\documentclass[twoside,draft,reqno,12pt]{amsart}
\usepackage[english]{babel}
\usepackage{amsmath}
\usepackage{amssymb}
\usepackage{enumerate}

\usepackage{amsfonts}
\usepackage{graphicx}

\usepackage[ citecolor=black, linkcolor=black,
colorlinks, hypertexnames=false, plainpages=false]{hyperref}

\usepackage{color}

\newtheorem{lemma}     {Lemma}[section]
\newtheorem{thm}   [lemma]{Theorem}
\newtheorem{teorema1}   [lemma]{Theorem}
\newtheorem{prop}       [lemma]{Proposition}

\newtheorem{cong1}      [lemma]{Conjecture}
\newtheorem{remark1}    [lemma]{Remark}

\numberwithin{equation}{section}

\newcommand{\und}{\underline}

\renewcommand{\(}{\left(}        \renewcommand{\)}{\right)}
\renewcommand{\[}{\left[}        \renewcommand{\]}{\right]}
     \newcommand{\nn}{\nonumber}

  \newcommand{\M}{{\mathcal M}}

\newcommand{\dis}{\displaystyle}

\newcommand{\mmmintone}[1]{{\dis{\int\kern -.38cm
-}}_{\kern-.21cm\substack{#1}}\;\;}
\newcommand{\mmmintwo}[2]{{\dis{\int\kern -.43cm
-}}_{\kern-.21cm\substack{#1}}^{\substack{#2}}\;\;}
\newcommand{\submint}{{\scriptstyle{\int\kern -.66em -}}}
\newcommand{\submintone}[1]{{\scriptstyle{\int\kern -.66em
-}}_{\scriptscriptstyle{\kern-.21em\substack{#1}}}}
\newcommand{\fracmint}{{\textstyle{\int\kern -.88em -}}}
\newcommand{\fracmintone}[1]{{\textstyle{\int\kern -.88em
-}}_{\scriptscriptstyle{\kern-.21em\substack{#1}}}\;}

\newcommand{\eps}{\epsilon}

\newcommand{\ga}{\gamma}
\newcommand{\Ga}{\Gamma}
\newcommand{\Om}{\Omega}
\newcommand{\om}{\omega}
\newcommand{\si}{\sigma}

\newcommand{\La}{\Lambda}

\newcommand{\E}{\mathbb E}
\newcommand{\Pp}{\mathbb P}
\newcommand{\PP}{\mathbf P}
\newcommand{\EE}{\mathbf E}

\newcommand{\nada}[1]{}

\begin{document}
\today

\vskip.5cm
\title
{Truncated correlations in the stirring process with births and deaths}

\author{A. De Masi}
\address{Anna De Masi,
Dipartimento di Matematica, Universit\`a di L'Aquila \newline
\indent L'Aquila, Italy}
\email{demasi@univaq.it}

\author{E. Presutti}
\address{Errico Presutti,
Dipartimento di Matematica, Universit\`a di Roma Tor Vergata \newline
\indent Roma, 00133, Italy}
\email{Presutti@mat.uniroma2.it}

\author{D. Tsagkarogiannis}
\address{Dimitrios Tsagkarogiannis,
Dipartimento di Matematica, Universit\`a di Roma Tor Vergata \newline
\indent Roma, 00133, Italy}
\email{tsagkaro@mat.uniroma2.it}

\author{M.E. Vares}
\address{ Maria Eulalia Vares,
Centro Brasileiro de Pesquisas Fisicas Rua Xavier Sigaud, 150.\newline
     \indent  Rio de Janeiro, RJ Brasil
     22290-180}
\email{eulalia@cbpf.br}
\begin{abstract}
We consider the stirring process in the interval $\La_N:=[-N,N]$ of $\mathbb Z$ with  births and deaths taking place in the intervals
$I_+:=(N-K,N]$, $K>0$, and respectively $I_-:=[-N,-N+K)$. We prove bounds on the truncated moments uniform in $N$ which yield strong factorization properties.

\end{abstract}

\maketitle

\tableofcontents

\section{Introduction}

We consider a system of particles in the interval $\La_N:=[-N,N]$ of
$\mathbb Z$. Particle configurations are elements $\eta$ of
$\{0,1\}^{\La_N}$, $\eta(x)=0,1$ being the occupation number at
$x\in\La_N$.  The evolution is a Markov process defined by the
generator
      \begin{equation}
        \label{pro.1.1}
L_\eps:=\eps^{-2}\big(L_0+\eps L_b\big), \quad \eps\equiv \frac 1N,
     \end{equation}
where $L_0$ is the generator of the stirring process in $\La_N$, namely
      \begin{equation}
        \label{pro.1.2}
 L_0 f(\eta):=\frac 12\sum_{x=-N}^{N-1} [f(\eta^{(x,x+1)})-f(\eta)]
     \end{equation}
with $\eta^{(x,x+1)}$ the configuration obtained from $\eta$ by interchanging
its values in $x$ and $x+1$, namely  $\eta^{(x,x+1)}(z)=\eta(z)$ if
$z\notin\{x,x+1\}$ and  $\eta^{(x,x+1)}(x)=\eta(x+1)$,
$\eta^{(x,x+1)}(x+1)=\eta(x)$. Moreover calling $\eta^{(x)}$ the configuration
obtained from $\eta$ by changing its value in $x$, namely $\eta^{(x)}(z)=\eta(z)$ if $z\ne x$ and
 $\eta^{(x)}(x)=0,1$ if  $\eta (x)=1,0$, then
      \begin{eqnarray}
        \label{pro.1.3}
&& L_b= L_{b,+}+ L_{b,-},\;\;L_{b,\pm} f(\eta):= \frac{j}{2}
\sum_{x\in I_\pm}D_{\pm}\eta(x) [f(\eta^{(x)})-f(\eta)\Big],\;\;j>0\\&&
I_+=[N-K+1,N],\; \; I_-=[-N,-N+K-1]\nn\\
&& D_+u (x)= [1-u(x)]u(x+1)u(x+2)\cdots u(N), \quad  x\in I_+
       \nn\\
&& D_-u (x)=  u(x)[1-u(x-1)][1-u(x-2)]\cdots[1-u(-N)], \quad  x\in I_-\nn
     \end{eqnarray}

$N^2L_0$ describes a process where nearest neighbor sites exchange
their content at rate $N^2/2$, the factor $N^2$ ensures that
information propagates to the whole $\La_N$ in time of order one,
with positive probability uniformly in $N$. $NL_{b,+}$ describes a
birth process: at rate $jN/2$ a particle is created in the first
(starting from $N$) empty site in $I_+$, if no site is empty the
birth is aborted. Symmetrically  $NL_{b,-}$ describes deaths: at
rate $jN/2$ the first particle (starting from $-N$) in $I_-$ is
removed, if $I_-$ is empty there is no death. Our results will show
that the factor $N$ is the correct one to match the stirring rate.

\vskip.5cm

This can be viewed as a model in queuing theory.  ``Files'' are sent in to the
first available ``server'' in $I_+$, they are elaborated and then finally sent out
from a server in $I_-$.  While the ``public'' only sees how many files enter and come
out, what really happens to the files is just that they go from one server to the other
back and forth (as described by the stirring process) in a random fashion.  This pessimistic
view of the action of the servers is unfortunately not too unrealistic as experienced by some
of the authors.

A more physical interpretation of the model is in terms of mass transport: mass is sent in from the right into $I_+$, it diffuses along $\La_N$ and then comes out from $I_-$.  In both cases the relevant questions are: how efficient is the system ? namely, once a steady regime has established,  what is the actual particle flux as a function of the external parameter $j$ which rules the birth-death mechanisms in $I_{\pm}$.  Moreover, how does the system relax from some initial state ?  As the system is a D\"oblin  chain there is a unique invariant measure and convergence is exponentially fast, but we look for bounds and estimates uniform in $N$  as $N\to \infty$.

\vskip.5cm

In this paper we shall study the evolution till  times which grow like $\tau \log N$, $\tau>0$ suitably small and shall prove strong factorization properties for the evolution starting from any single configuration, (we prove rather sharp bounds on the truncated correlation functions). This result together with those  in a companion paper, \cite{DPTVjsp},
prove convergence of empirical averages to a limit macroscopic equation, see the next section for details. Our project for the future is to prove (i) that the time-flow defined by the limit macroscopic equation leads to a unique
stationary profile as $t\to \infty$; (ii) that the unique stationary measure of the process when $N$ is finite is in the limit $N\to \infty$ supported by the   macroscopic stationary profile. The strong factorization properties proved here seem to be the main ingredients for this program to be fulfilled, but the analysis of the limit equation has still to be done.

\vskip2cm

\section{Main result}

We shall study the process described in the introduction starting from an arbitrary
initial configuration $\eta$, denoting by $\Pp_\eps$ its law and $\E_\eps$ the corresponding expectation. We shall not make explicit the dependence on $\eta$ unless ambiguities may arise.  We are interested in the expectations
$\E_\eps[\eta(x,t)]$ and in the truncated correlation functions $\dis{\E_\eps\Big[\prod_{i=1}^n \tilde\eta(x_i,t)\Big]}$, $ \tilde\eta(x,t):= \eta(x,t)-\E_\eps[\eta(x,t)]$ (it is actually more convenient to study a slightly different expression as the one defined in \eqref{pro.2.3.1} below).  Define first $\rho_\eps(x,t)$ as the solution of
            \begin{eqnarray}
          \label{pro.2.1}
&&\frac{d}{dt}\rho_\eps(x,t)= \frac 12 \Delta_\eps\rho_\eps(x,t)+  \eps^{-1} \frac j2\Big(\mathbf 1_{x\in I_+} D_+\rho_\eps(x,t)
-\mathbf 1_{x\in I_-}  D_-\rho_\eps(x,t)\Big)\\&&
\rho_\eps(x,0)=\mu_\eps[\eta(x,0)=1]\nn
         \end{eqnarray}
where $\Delta_\eps=\eps^{-2}\Delta$, $\Delta$
the discrete Laplacian in $\La_N$ with reflecting boundary conditions:
    \begin{eqnarray}
   \nn
&&\Delta u(x)=   u(x+1)+ u(x-1)-2u(x),\qquad |x|<\eps^{-1}
\\&&\Delta u(\pm N)= u(\pm(N-1),t)-u(\pm N,t)
         \label{pro.2.2}
        \end{eqnarray}
 and $\mu^\eps$ is a product measure [in particular it may be supported by a single configuration].
Global existence and uniqueness for \eqref{pro.2.1} are proved in \cite{DPTVjsp} where it is also shown that the solution has values in $[0,1]$.  Writing $\La_N^{n,\ne}$, $n\ge 1$, for the set of all sequences
 $\und x=(x_1,..,x_n)$ in $\La_N^n$ such that
 $x_i\ne x_j$ we then define the $v$-functions
             \begin{equation}
             \label{pro.2.3.1}
v^\eps(\und x,t|\mu^\eps):= \E_\eps\Big[\prod_{i=1}^n \{ \eta(x_i,t)-\rho_\eps(x_i,t)\}
\Big], \quad \und x\in \La_N^{n,\ne},\; n\ge 1
             \end{equation}
where the process starts with a product measure $\mu^\eps$
and $\rho_\eps(x,t)$ is the solution of \eqref{pro.2.1}
The main result in this paper is a bound on the $v$-functions:

\vskip.5cm

                \begin{thm}
                \label{pro.thm.2.3}

There exist $\tau>0$ and $c^*>0$ so that the following holds. For any $\beta^*>0$
and  for any positive integer $n$ there is a constant $c_n<\infty$ so that for any
$\eps>0$, any initial product measure $\mu^\eps$
             \begin{equation}
             \label{pro.2.3}
\sup_{\und x\in \La_N^{n,\ne}}|v^\eps(\und x,t|\mu^\eps)| \leq \begin{cases}
c_n(\eps^{-2}t)^{-c^* n}, & t\le \eps^{\beta^*}\\
 c_n \eps^{(2-\beta^*)c^* n} &   \eps^{\beta^*} \le t\le \tau\log\eps^{-1}.
 \end{cases}
             \end{equation}
             \end{thm}

\vskip.5cm
Theorem \ref{pro.thm.2.3} proves that
$\sup_x|\E_\eps[\eta(x,t)]-\rho_\eps(x,t)| \to 0$ as $N\to \infty$
and that the ``empirical averages'' behave as $\rho_\eps$ in the
following sense. Call $J_N(x)=[x-N^a,x+N^a]\cap \La_N$,
for $a\in (0,1)$, and let $|A|$ denote the cardinality of $A \subset \mathbb Z$.
Then by the Chebyshev inequality,
             \begin{equation}
             \label{pro.2.4}
\lim_{\delta\to 0}\lim_{\eps\to 0}\sup_{\eps^{\beta^*}\le t \le \tau
\log\eps^{-1}} \sup_{\eta} \Pp_\eps\Big[\sup_{x\in \La_N}
\big|\frac{1}{|J_M(x)|}\sum_{y\in J_M(x)}\{\eta(y,t)-\rho_\eps(y,t)\}\big|\ge
\delta
 \Big]=0.
              \end{equation}
Then under suitable assumptions on the initial configurations it follows that the above empirical averages converge to the solution of the hydrodynamic equation for the system, as we are going to explain.  Suppose that the initial configurations $\rho_\eps(x,0)$ (see the second line in \eqref{pro.2.1}) are such that for some smooth function $u_0(r), r\in [-1,1]$,
             \begin{equation}
             \label{pro.2.5}
\lim_{\eps\to 0} \sup_{x\in \La_N}  \Big|\frac1{|J_M(x)|} \sum_{y\in J_M(x)}\{\rho_\eps(y,0)-u_0(\eps y)\}\Big|=0.
             \end{equation}
In \cite{DPTVjsp} it has been proved that under such an assumption for any $t>0$,
             \begin{equation}
             \label{pro.2.6}
\lim_{\eps\to 0} \sup_{x\in \La_N} | \rho_\eps(x,t)-\rho(\eps x,t) | =0
             \end{equation}
with $\rho(r,t)$ the unique solution of the limit hydrodynamic equation, namely the heat equation  with Dirichlet boundary conditions:
             \begin{eqnarray}
             \label{pro.2.7}
&&\frac {\partial}{\partial t} \rho(r,t)= \frac 12 \frac {\partial^2}{\partial r^2} \rho(r,t), \qquad r\in (-1,1), t>0 \\&& \rho(r,0)=u_0(r),\quad \rho(\pm 1,t)= u_{\pm}(t)\,.\nn
             \end{eqnarray}
However the boundary conditions $u_{\pm}(t)$ are not a-priori known, they must be obtained by solving a nonlinear system of two integral equations:
             \begin{eqnarray}
             \label{pro.2.8}
&& u_{\pm}(t) = \int_0^t\{p(s) f_{\pm}(u_{\pm}(t-s))-q(s) f_{\mp}(u_{\mp}(t-s)) \} ds + w_{\pm,t}\\&&
f_+(u)= \frac j2\Big(1- u^K\Big),\quad f_-(u)= \frac j2\Big(1 - (1-u)^K\Big)\,, \nn
             \end{eqnarray}
 where
              \begin{eqnarray}
             \label{pro.2.9}
&& p(t) = 2\sum_{k\in \mathbb Z} G_t(4k),\;\;q(t) = 2\sum_{k\in \mathbb Z} G_t(4k+2),\;\;\; G_t(r)= \frac{e^{-r^2/2t}}{\sqrt{2\pi t}}\nn\\
\\&&
w_{+,t} = \sum_{k\in \mathbb Z} \int_{-1}^1 u_0(r') 2 G_t(1-r' +4k)dr',\;\;\; w_{-,t} = \sum_{k\in \mathbb Z} \int_{-1}^1 u_0(r') 2 G_t(r'+1 +4k)dr'.\nn
             \end{eqnarray}
Theorem \ref{pro.thm.2.3} then shows that \eqref{pro.2.4} holds with $\rho_\eps( y,t)$ replaced by $\rho(\eps y,t)$ and for $t$ in a compact interval.

\vskip.5cm
{\bf Scheme of proofs}. In Section \ref{pro.sec.3} we shall state some mostly elementary properties of simple random walks on $\La_N$ which have been proved or recalled in \cite{DPTVjsp}. In Section \ref{secN6} we shall prove sharp probability estimates on the stirring process, which extend analogous estimates proved in \cite {DP} for the process on the whole $\mathbb Z$. The reflections at $\pm N$ make the extension not trivial at all.   In the remaining sections we write an integral equation for some truncated correlation functions called ``the $v$-functions'' that we study by iteration. For further results and applications of $v$-functions see e.g. \cite {BPSV,DPS,FPSV}.
The terms which arise are interpreted as a branching process with stirring evolution between the branching events.  The main point will be to prove that typically branching events are well separated in time so that the stirring has time
to mix up things in the proper way, here we use extensively the estimates  in Section  \ref{secN6}.

\vskip2cm

\section{A single random walk with reflections}
\label{pro.sec.3}

In this section we state some properties of a single random walk in $\La_N$ with reflections
at $\pm N$, referring to \cite{DPTVjsp} for the proofs.

We denote by
$P_t^{(\eps)}(x,y)$   the transition probability of a simple random walk in $\La_N$
which jumps with intensity $\eps^{-2}/2$ to each of its n.n. sites,
jumps outside  $\La_N$ are suppressed. $Q_t^{(\eps)}(x,y)$ is instead the transition
probability of a random walk on the whole $\mathbb Z$.
$P_t^{(\eps)}$ and  $Q_t^{(\eps)}$ are related via the  ``reflection map'' $\psi_N:\mathbb Z \to \La_N$ defined as follows.
\vskip.5cm

\begin{itemize}

 \item\; $|x|\le N$:\;\;$\psi_N(x)=x$.

 \item\; $x<-N$:\;\;$\psi_N(x)=-\psi_N(-x)$;

 \item\; $x>N$:\;\;$\psi_N\left(N+j(2N+1)+k\right)=\psi_N\left(N+j(2N+1)-(k-1)\right)$,\;
  $k=1,\dots,2N+1, j=0,1,\dots$.

\end{itemize}

\vskip.5cm

Then
     \begin{equation}
    \label{a3.6}
P^{(\eps)}_t (x,z)= \sum_{y: \psi_N(y)=z} Q^{(\eps)}_t (x,y).
    \end{equation}
   Let
      \begin{equation}
G_t (x,y)=   \frac{e^{- (x-y)^2/2t}}{\sqrt{2\pi t}},
    \end{equation}
By the local central limit theorem, \cite{lawler},
there exist positive constants $c_1,...,c_5$ so that
    \begin{equation}
    \label{N4.4}
 |Q^{(\eps)}_t(x,y)- G_{\eps^{-2}t}(x,y)|\le   \frac{c_1}{ \sqrt{\eps^{-2}t}} G_{\eps^{-2}t}(x,y),\quad |x-y|\le (\eps^{-2 }t)^{5/8}.
    \end{equation}
 while for $|x-y|> (\eps^{-2 }t)^{5/8} $
    \begin{equation}
       \label{N4.4.1}
 Q^{(\eps)}_t(x,y) \le \min\Big\{ c_2 e^{-c_3|x-y|^{2}/(\eps^{-2}t)}, c_4 e^{-|y-x|(\log |y-x|-c_5)}\Big\}.
    \end  {equation}

\vskip1cm

As a corollary, for any $T>0$ there exists $c$  so that the following holds:

\begin{itemize}

\item For all $\eps$, all $t\in (0,T]$ and all $x$, $y$ in $\La_N$
     \begin{equation}
    \label{N4.5}
 P^{(\eps)}_t (x,y)    \le c \;G_{\eps^{-2}t}(x,y).
    \end{equation}
\item For all $\eps$, all $t\in (0,T]$ and all
 $-N\le x \le N-1$,
     \begin{equation}
    \label{N4.6}
\Big| P^{(\eps)}_t (x,y) - P^{(\eps)}_t (x+1,y) \Big| \le \frac{c}{\sqrt{\eps^{-2}t}} G_{\eps^{-2}t}(x,y).
    \end{equation}

    \end{itemize}

\vskip.5cm
In Proposition 5.1 of \cite{DPTVjsp} it is proved that
for any $S>0$ there is a constant $c$ so that for any solution $\rho_\eps(x,t)$ of \eqref{pro.2.1} with $\rho_\eps(\cdot,0)\in [0,1]$
the following holds. For any $x\in[-N,N-1]$, any $t\in (0,S]$  and any $\eps>0$
       \begin{equation*}
       \label{N5.2}
|\rho_\eps(x,t)-\rho_\eps(x+1,t)| \le \min\Big\{ 1, c  \Big( \eps  \log_+ (\eps^{-2}t) + \frac{1}{\sqrt{\eps^{-2}t}}\Big)\Big\}
                \end{equation*}
where $ \log_+ u =\max\{\log u, 1\} $.  By the arbitrariness of the initial datum we can iterate
the bound and  from \eqref{N5.2} with $S=1$ we  get that for all $x\in[-N,N-1]$ and all $t>0$:
       \begin{equation}
              \label{N5.2a}
|\rho_\eps(x,t)-\rho_\eps(x+1,t)| \le \min\Big\{ 1, c  \Big( \eps  \log_+ (\eps^{-2}[t]_1) + \frac{1}{\sqrt{\eps^{-2}[t]_1}}\Big)\Big\}
                \end{equation}
where $[t]_1=t$ if $t\le 1$ and $=1$ otherwise.  We shall use a weaker version of \eqref{N5.2a}, namely that
  for any $\zeta>0$ and $\tau>0$ there is $c$ so that
       \begin{equation}
       \label{18a}
 \sup_{x,y\in\La_N: |x-y|\le 1}|\rho_\eps(x,t)-\rho_\eps(y,t)| \le \frac {c}{(\eps^{-2}t)^{1/2-\zeta}+1}
\quad \text{{\rm for any} $t\le \tau \log \eps^{-1}$}.
                \end{equation}

%
%
%

\vskip2cm

\section{Probability estimates for the stirring process.}
    \label{secN6}

In this section we shall study the stirring process.
The generator $L_0$ of the stirring process defined in \eqref{pro.1.2} can be interpreted by saying that after independent exponential times of mean $1/2$ sites $x$ and $x+1$ exchange their content. This leads to the following realization of the process which will be extensively used in the sequel.

\vskip.5cm

\noindent
{\bf Definition 1.} [{\em   The active/passive marks process}]

$\bullet$\; The active/passive marks process is realized in a probability space denoted by $(\Om,\PP_\eps)$. It is a product of Poisson processes indexed by $\{x,x+1\}$, $x\in \mathbb Z$: for each pair $\{x,x+1\}$ we have a Poisson point process of intensity $\eps^{-2}$, its events are called ``marks'' and each mark is independently given the attribute ``active'' or ``passive'' with probability $1/2$. The processes relative to different pairs are mutually independent and their common law is $\PP_\eps$, its expectation being also denoted by $\EE_\eps$.

$\bullet$\;  For any $\om\in \Om$ we  define the following particle evolution in $\La_N$: a particle at $x\in \La_N$ moves as soon as an active mark appears at a pair $x,y$ with $y\in \La_N$. The particle then moves to $y$ (note that if another particle was at $y$ then it would simultaneously move to $x$). Passive marks do not play any role so far as well as all marks  $\{x,x+1\}$ with at least one site not in $\La_N$, they will be used later to construct couplings.  We denote by $X(t)\subset \La_N$ the set of all sites occupied by the  particles  at time $t$, so that $\eta(x,t)=1$ iff $x\in X(t)$.

\vskip.5cm

\begin{prop}
\label{propN6.1}

$X(t)$ has the law of the stirring process with generator $\eps^{-2}L_0$.

\end{prop}

\vskip.5cm
The  proof of the proposition is elementary and omitted.  The evolution defined in the active/passive marks process can clearly be inverted: given $\om$ and $X(t)$ by following backwards the marks we uniquely determine the initial position $X(0)$. This remark leads to the proof of the well known proposition, see for instance \cite{liggett}:

\vskip.5cm

\begin{prop} [Duality]
\label{propN6.2}
For any $X\subset \La_N$, $\eta_0\in\{0,1\}^{\La_N}$ and $t>0$,
              \begin{equation}
              \label{N6.2}
 \EE_\eps\Big[ \prod _{x\in X} \eta(x,t)\;|\; \eta(\cdot,0)=\eta_0\Big] = \EE_\eps\Big[ \prod _{x\in X(t)} \eta_0(x)\;|\; X(0)=X\Big]
              \end{equation}

\end{prop}

\vskip.5cm

The importance of Proposition \ref{propN6.2} is that it allows to compute the probability at time $t$ that the sites in $X$ are all occupied by studying the
stirring process with ``only'' $|X|$ particles (no matter how many are the particles in the initial configuration $\eta_0$).  The realization of the process in terms of   active/passive marks allows to identify particles and hence to label them:

\vskip.5cm

\noindent
{\bf Definition 2.} [{\em   The labeled process}]

  Given $\om\in \Om$, we can follow unambiguously the motion of each individual particle so that we can give them   labels at time 0  which then remain attached to the particles during their motion.
  We shall denote by $\und x=(x_{i_1},\dots,x_{i_n})$ a labeled configuration of $n$ particles, ($i_1,\dots,i_n$ the labels, $x_{i_j}$ the positions); configurations obtained under permutations of the labels are now considered  distinct. We write $\und x(t)$
for the labeled process induced by the active/passive marks and denote by $X(t)$ the unlabeled configuration obtained from $\und x(t)$, (i.e.\ only the positions in $\und x(t)$ are recorded by $X(t)$).
By an abuse of notation we shall write $\PP_\eps$ and $\EE_\eps$ both for law and expectation in the active/passive marks process and for the marginal over   $\und x(t)$, the labeled process realized in this space.  By adding a subscript $\und x$ we mean that  the initial distribution of particles has support on the single labeled configuration $\und x$.

\vskip.5cm

The advantage of having defined the process with also passive marks is exemplified in Lemma \ref{x}  below, where we use the following:

\vskip.5cm

\noindent
{\bf Definition 3.} [{\em   The variables $\mathcal T_{x_1,x_2}$, $\tau_{x_1,x_2}$ and  $N_{x_1,x_2,t}$}]

Given the initial position  $x_1$ and $x_2$ of particles $1$ and $2$ define in $\Om$ the random multi-interval $\mathcal T_{x_1,x_2}=\{s\ge 0: |x_1(s)-x_2(s)|=1\}$,  calling $\mathcal T_{x_1,x_2,t}:=\mathcal T_{x_1,x_2}\cap [0,t]$ and $I_{x_1,x_2,t}=\{ (s,y_1,y_2): s \in \mathcal T_{x_1,x_2,t}, y_i=x_i(s)$ and a mark appears at $s$ between $y_1$ and $y_2$$\} $,
we define the stopping time $\tau_{x_1,x_2}$ as the smallest
$s$ in $I_{x_1,x_2,\infty}$ 
 and    $N_{x_1,x_2,t}$   the total number of elements in $I_{x_1,x_2,t}$.
Notice that the presence of other particles does not affect the values  of $\mathcal T_{x_1,x_2}$, $\tau_{x_1,x_2}$ and  $N_{x_1,x_2,t}$.

\vskip.5cm

      \begin{lemma}
      \label{x}
Let $\und x=(x_1,\dots,x_n)$, $t>s>0$ and $f(y_1,\dots,y_n)$ a function antisymmetric under the exchange of $y_1$ and $y_2$.
Then
       \begin{equation}
       \label{18a.1}
 \EE_{\eps,\und x}\Big[ \mathbf 1_{\tau_{x_1,x_2} \le s } f(\und x(t))\Big] = 0,
                \end{equation}
where the suffix $\und x$ indicates the initial condition i.e. $\EE_{\eps,\und x}(\cdot)=\EE_{\eps} [\cdot|\und x(0=\und x]$.

      \end{lemma}

\vskip.5cm

 {\em Proof.} Given $x_1$,  $x_2$ and $s>0$ call
two  elements $\om$ and $\om'$ in $\Om$   ``similar'' if all the marks at all pairs $x,x+1$ occur at the same time in $\om$ and $\om'$; their active/passive attribute must also be the same in both except for the marks in $I_{x_1,x_2,s}$  (see Definition 3): a mark at the times and between the pairs indicated by  $I_{x_1,x_2,s}$ may be active in one sample and passive in the other, or the same in both.
If $\om$ and $\om'$ are similar then the configurations evolved from the same $\und x$ with $\om$ and $\om'$ are  at all times the same except at most for an exchange of the particles with label 1 and 2.  Then the above similarity relation is an equivalence. $N_{x_1,x_2,s}$ is constant in all $\om$ in an equivalence class, so that if $N_{x_1,x_2,s}=p$ there are $2^p$ elements
in the equivalence class.  Each element is characterized by the active/passive attribute of the $p$ marks  in $I_{x_1,x_2,s}$ and their distribution,  conditioned on a given equivalence class,
is a product of  $\frac 12, \frac 12$ probabilities.
  The set $\{\tau_{x_1,x_2} \le s\} $ is  the union of all such equivalence classes with $p\ge 1$ so that
  the law  of $\und x(t)$ conditioned on $\{\tau_{x_1,x_2} \le s\} $ is symmetric under exchanging $x_1(t)$ with $x_2(t)$ and \eqref{18a.1} follows.

\qed

\vskip1cm

\noindent

We shall complement Lemma \ref{x} by proving in \eqref{19.e1} below
that the probability that $\tau_{x_1,x_2} \le s$ goes to 1 as $s\to \infty$ if $|x_1-x_2|=1$, thus establishing
exchangeability properties of the stirring process.
Since the stirring particles move like independent random walks when at distance larger than 1, the crucial ingredients to construct good couplings between the two processes (as in \cite{FPSV, DP}) are a-priori estimates on the tails of the variables $\tau_{x_1,x_2}$ and $N_{x_1,x_2,t}$.  The results of \cite{DP} cannot be directly applied
here because they heavily used that the particles were moving on the whole line, and our process is in the bounded interval $\La_N$. Some estimates are helped by being in a bounded domain but in others the inequality goes in the wrong direction and we have extra work to do.

 \vskip1cm

   \begin{thm}
   \label{3.01}
There is  $c$ so that for all $\eps>0$
                    \begin{equation}
              \label{19.e1}
\sup_{|x_1-x_2|=1} \PP_\eps\big[\tau_{x_1,x_2} \ge t\big] \le \frac{c}{(\eps^{-2}t)^{1/2}+1}.
             \end{equation}
Moreover given any $T>0$ for any $\zeta>0$ and $k\ge 1$ there is $c$ so that for all $t\le T$ and  for all $\eps>0$
                      \begin{equation}
              \label{19.e2}
\sup_{x_1,x_2} \PP_\eps\big[N_{x_1,x_2,t} \ge (\eps^{-2}t)^{1/2+\zeta}\big] \le  c (\eps^{-2}t)^{-k}.
             \end{equation}

\end{thm}

\vskip.5cm

{\em Proof.} The proof of \eqref{19.e1} will follow by bounding $\PP_\eps\big[\tau_{x_1,x_2} \ge t\big]$ in terms of the probability of the return time to the origin of a single random walk on $\mathbb Z$.
Given a realization $\om$ of the active/passive marks process, call $x_i(t)$, $i=1,2$, the positions of the particles evolved from $x_1$ and $x_2$ by considering only the marks in $\La_N$ (as in Definition 3 above); denote instead by $x^+_i(t)$, $i=1,2$, the positions obtained by using all the marks in $\mathbb Z$ (thus $x^+_i(t)$ are two stirring walks on the whole $\mathbb Z$ starting from $x_1$ and $x_2$). Suppose without loss of generality that $x_1<x_2$ then for all $t\le \tau_{x_1,x_2}$ we have $x^+_1(t)\le x_1(t) < x_2(t)\le x^+_2(t)$ so that 
$\PP_\eps[\tau_{x_1,x_2}\ge t]$ is bounded from above by the probability that $t$ is smaller than the first time $s$ when two independent random walks on $\mathbb Z$ starting from $x_1$ and $x_2$ are on a same site. Thus \eqref{19.e1} follows from classical estimates on return times of random walks.

The l.h.s.\ of \eqref{19.e2} is bounded by
           \begin{equation*}
\EE_\eps\Big[ \PP_\eps\Big[N_{x_1,x_2,t} \ge (\eps^{-2}t)^{1/2+\zeta}\;\big|\;|\mathcal T_{x_1,x_2,t}| < \eps^{2}(\eps^{-2}t)^{1/2+\zeta'} \Big]\,\Big] + \PP_\eps\Big[|\mathcal T_{x_1,x_2,t}| \ge \eps^{2}(\eps^{-2}t)^{1/2+\zeta'} \Big].
             \end{equation*}
Let   $\dis{p_L(n)= e^{-L} \frac{L^n}{n!}}$ the Poisson distribution  with parameter $L$, then
                        \begin{equation*}
 \PP_\eps\Big[N_{x_1,x_2,t} \ge (\eps^{-2}t)^{1/2+\zeta}\;\big|\;|\mathcal T_{x_1,x_2,t}| =L \Big] = \sum_{n\ge  (\eps^{-2}t)^{1/2+\zeta}} p_{\eps^{-2}L}(n)
             \end{equation*}
Since the r.h.s.\ is an increasing function of $L$,
             \begin{eqnarray*}
&&\hskip-1cm\EE_\eps\Big[ \PP_\eps\Big[N_{x_1,x_2,t} \ge (\eps^{-2}t)^{1/2+\zeta}\;\big|\;|\mathcal T_{x_1,x_2,t}| < \eps^{2}(\eps^{-2}t)^{1/2+\zeta'} \Big]\,\Big]\\&&\hskip4cm  \le \sum_{n\ge(\eps^{-2}t)^{1/2+\zeta}} p_{(\eps^{-2}t)^{1/2+\zeta'}} ( n)
             \end{eqnarray*}
which, for any $k$, is at most $c(\eps^{-2}t)^{-k}$ if  $\zeta'=\zeta/2$ and    $c$ is a suitable constant (dependent on $k$ and $\zeta$).
It will therefore suffice to prove that for any $\zeta>0$ and any $k$ there is $c$ so that
  \begin{equation}
              \label{19.e3}
\sup_{x_1,x_2} \PP_{\eps,x_1,x_2}\big[\int_0^t \mathbf 1_{|x_1(s)-x_2(s)|=1}ds \ge \eps^{2}(\eps^{-2}t)^{1/2+\zeta}\big] \le  c (\eps^{-2}t)^{-k},
             \end{equation}
where as before the suffix $x_1,x_2$ in $\PP_{\eps, x_1,x_2}$  recalls that $x_i(0)=x_i, i=1,2$. By the Chebyshev inequality:
                  \begin{equation}
              \label{19.e4}
\text{l.h.s. of \eqref{19.e3}}    \le \sup_{x_1,x_2}\big(\eps^{-2}t\big)^{-(\frac 12 +\zeta)p} p!\int_{\{0\le t_1 \le t_2\dots \le t_p \le t\}} \eps^{-2}dt_1...\eps^{-2}dt_p
\EE_{\eps,x_1,x_2} \Big[\prod_{i=1}^p\mathbf 1_{|x_1(t_i)-x_2(t_i)|=1}\Big].
             \end{equation}
By choosing $p$ so that $\zeta p>k$ we are left with the proof that there is a constant $c$ so that
                  \begin{equation}
              \label{19.e5}
\sup_{x_1,x_2}\int_{\{0\le t_1 \le t_2\dots \le t_p \le t\}} \eps^{-2}dt_1...\eps^{-2}dt_p
\EE_{\eps,x_1,x_2} \Big[\prod_{i=1}^p\mathbf 1_{|x_1(t_i)-x_2(t_i)|=1}\Big] \le c (\eps^{-2}t)^{p/2}.
             \end{equation}
By conditioning successively on the state at the times $t_{p-1}, t_{p-2},\dots,0$ we are reduced to prove the bound
                  \begin{equation}
              \label{19.e6}
\sup_{y_1,y_2}\PP_{\eps,y_1,y_2} \Big[{|x_1(s)-x_2(s)|=1}\Big] \le \frac c {(\eps^{-2}s)^{1/2}}
             \end{equation}
where $s$ stands for $t_{i+1}-t_i$.  \eqref{19.e5} follows from  \eqref{19.e6} after some simple computations (details are omitted) and we are left with the proof of \eqref{19.e6}.
We use the Liggett's inequality to bound
                 \begin{equation}
              \label{19.e7}
\PP_{\eps,y_1,y_2}\Big[{|x_1(s)-x_2(s)|=1}\Big] \le \sum_{x} \PP_{\eps,y_1}\Big[x_1(s)\in\{x,x+1\}\Big] \PP_{\eps,y_2}\Big[x_2(s)\in\{x,x+1\}\Big].
             \end{equation}
\eqref{19.e6} then follows using \eqref{N4.4} and \eqref{N4.4.1}. \qed

\vskip.5cm

The above estimates will be used to prove that $n$ stirring particles move like $n$ independent random walks.  We
construct a coupling between the two processes by realizing both of them in the active/passive mark processes.  The definition
adapts to the present case the one considered  in \cite{DP}.

\vskip1cm

{\bf Definition 4.} [Couplings]

Denoting by  $\und x$ and $\und x^0$  stirring and independent labeled particles, respectively, without loss of generality we suppose that the  labels are  $1,\dots,n$  and write $\und x= (x_1,\dots,x_n)$, $\und x^0= (x_1^0,\dots,x_n^0)$. We also assign (arbitrarily) a priority list $\si$, $\si$ a permutation of $\{1,\dots,n\}$, and say that particle $i$ has priority over particle $j$ if $\si(i)<\si(j)$. We consider the
active/passive marks process realized on the whole $\mathbb Z$ and according to Definition 1 in this Section we define $\und x(t)$ by looking only at the marks $\{x,x+1\}$ with both $x$ and $x+1$ in $\La_N$.
In the same space we define $\und x^0(t)= (x_1^0(t),\dots,x^0_n(t))\in \La_N^n$ given $\und x(0)$, $\und x^0(0)$ and $\si$ by the following rules (we shall later prove that this is a process of independent random walks in $\La_N$).
$\und x^0(t)$ is defined  by giving for each $i\in\{1,\dots,n\}$ the times $t_{i,r}$, $t_{i,l}$ when  $x^0_i(\cdot)$ ``tries'' to move to the right, respectively to the left, as those jumps which lead out of $\La_N$ are suppressed. Thus the sequence of all such times determines $\und x^0(\cdot)$ (while the viceversa is not true
as we cannot recover from $x^0(\cdot)$ the attempts to jump out of $\La_N$). While the times  $t_{i,r}$, $t_{i,l}$ cannot be recovered from $x^0$ they can be read out of   an auxiliary   process $\und y(t)$ that we define next.

 \vskip.5cm

{\em Definition of $y(t)$.} Let $t > 0$ be the first time when a mark appears at a pair $(x,x+1)$ such that $\und x(0)\cap\{x,x+1\}\ne \emptyset$. We set $\und y(s)=\und x^0(0)$ for $s\in (0,t)$ and  define $\und y(t)$ (which will then define  $x^0(s)$, $s\le t$, as well) in the following way.

\vskip.5cm

 CASE 1, $\und x(0)\cap \{N,-N\}=\emptyset$.

$\bullet$\;Subcase (1.i): the intersection $\und x(0)\cap\{x,x+1\}$
is a singleton, for instance $\und x(0)\cap\{x,x+1\}=x_i(0)$.

(1.i.a):  the mark is passive, then
$\und y(t)=\und x^0(0)$.

(1.i.b): the mark is active,  then  $y_k(t)=  x^0_k(0)$, $k\ne i$, and
$y_i(t)= x^0_i(0)\pm 1$ if $x_i(0)=x$, respectively if  $x_i(0)=x+1$.

$\bullet$\; Subcase (1.ii): the intersection is a doubleton, for instance $\und x(0)\cap\{x,x+1\}=\{x_i(0), x_j(0)\}$ and $\si(i)<\si(j)$; we then say that ``particles $i$ and $j$ collide at time $t$''. 

(1.ii.a): the mark is passive, then  $y_k(t)=  x^0_k(0)$, $k\ne j$ and  $y_j(t)=  x^0_j(0)-[x_j(0)-x_i(0)]$

(1.ii.b): the mark is active. Then   $y_k(t)=  x^0_k(0)$, $k\ne i$ and $y_i(t)=  x^0_i(0)+[x_j(0)-x_i(0)]$

\vskip.5cm

 CASE 2, $\und x(0)\cap \{N,-N\}\ne \emptyset$.


If both $x$ and $x+1$ are in $\La_N$   the same rules given before apply. It remains to consider the subcases where $x=-N-1$ and $x=N$, the definitions are analogous and we only consider the former case.  Let $x_i(0)=-N$, so that $x=-N-1$ and $x+1=x_i(0)$.

$\bullet$\;Subcase (2.i): The mark is passive, then $\und y(t)=\und x^0(0)$.

$\bullet$\;Subcase (2.ii): The mark is active, then  $y_k(t)=  x^0_k(0)$, $k\ne i$,
 $ y_i(t)= x^0_i(0)-1$. 

\vskip1cm
Having defined   $\und y(s)$ till $s=t$ we
then set $\und x^0(s)=\und y(s)= \und x^0(0)$, $s< t$; and $\und x^0(t)=\und y(t)$ if $\und y(t)\in\La_N^n$, otherwise  $\und x^0(t)=\und x^0(0)$.  Having  $\und x(t)$ and $\und x^0(t)$ we can then define $\und y(s)$  from $t$ till the time of the next mark by using the same rules used starting from time 0. In particular they imply that $\und y(t^+)=\und x^0(t)$ (so that it may happen that at time $t$
$\und y(\cdot)$ jumps twice if the first jump takes $\und y$ out of $\La_N$, then ``instantaneously'' $\und y$ comes back from where it jumped and have $\und y(t^-)=\und y(t^+)\ne \und y(t)$).  By applying repeatedly this procedure we define (with probability 1) $\und y(\cdot)$ and  $\und x^0(\cdot)$ at all times.  The collections of times when  $y_i$ jumps to its right/left are called respectively $\{t_{i,r}\}$ and $\{t_{i,l}\}$ (right/left refers to the first jump if there are two jumps at the same time).  These are the ``attempted jumps'' of $x^0_i$ because when  $y_i$ jumps twice (the first time out of $\La_N$ and the second time back to the initial position) then $x^0_i$ does not change, i.e.\ the jump is  suppressed.

\vskip1cm
It readily follows from the definition (see \cite{DP} for details)

\begin{itemize}


\item  The jump times $\{t_{i,r},t_{i,l}, i=1,\dots,n\}$  are mutually independent Poisson processes with intensity $\eps^{-2}/2$ and
 $\und x^0(t)$ has the same law as $n$ independent random walks in $\La_N$, with jump rate $\eps^{-2}/2$ to each n.n.\ site in $\La_N$.

\item  The times when $x_i$ and $x^0_i$ have different jumps can only occur when one of them is at $\pm N$ or when ``$x_i$ collides with $x_j$'' and $\si(j)<\si(i)$ (see subcase (1.ii)).

\item For any $i$ and $t>0$, $x_i(t)$ is completely determined by $y_j(s)$, $s\in [0,t]$ with $j:\si(j)\le \si(i)\}$.

\item  If $\si(\ell)=1$ ($\si(\cdot)$ the ``priority list'')
and   $x_{\ell}(0)=x_{\ell}^0(0)$ then, with probability 1, $x_{\ell}(t)=x_{\ell}^0(t)$ for all $t\ge 0$.

\end{itemize}

\vskip1cm
The following theorem proves bounds on $|x_i(t)-x_i^0(t)|$  as those established  in \cite{DP} for processes on $\mathbb Z$, but the proof is more involved due to the reflections at $\pm N$.

 \vskip.5cm

   \begin{thm}
   \label{5e.3}
Let $T>0$ and $\und x(0)= \und x^0(0)$. Then for any $\zeta>0$ and $k$ there is $c$ so that for all $t\le T$ and  for all $\eps>0$
                      \begin{equation}
              \label{19.e11}
\PP_\eps\big[|x_\ell(t)-x^0_\ell(t)|  \ge (\eps^{-2}t)^{1/4+\zeta}\big] \le  c (\eps^{-2}t)^{-k}
             \end{equation}

\end{thm}

\vskip.5cm

{\em Proof.}  Without loss of generality we may and will suppose that $\si$ is the identity permutation (just rename the particles following the priorities).  Since $x_1(t)=x_1^0(t)$ for all $t\ge 0$ we only need to prove \eqref{19.e11}
with $\ell>1$, the label $\ell$ being hereafter fixed.

 The idea of the proof is the following. Define the vectors
 $$
 D(t)=\big(x_1(t)-x^0_1(t),..,x_\ell(t)-x^0_\ell(t)\big), \quad
 \xi(t)=\big(x_1(t),..,x_\ell(t)\big)
 $$
 The first point  will be to check that $\dis{D^2(t) - \frac{\eps^{-2}}2H(\xi(t))}$ is a $\PP_\eps$ super-martingale, where $H(\xi)$ is the number of pairs $x_i,x_j$ in $\xi$ such that $|x_i-x_j|=1$.  We shall then need to extend the analysis to the higher moments  $D^{2n}(t)$ and in this way we shall relate the moments of $x_\ell(t)-x^0_\ell(t)$ to bounds on the probability  of the time length when pairs of stirring  particles are close-by.  All that however is more easily accomplished by studying the skeleton of the process, i.e.\ by looking at the times when particles jump and since the process is realized in the active/passive marks process, at the times when the marks appear.

Thus, following the proof of Lemma \ref{x}, we call two elements $\om$ and $\om'$  of $\Om$ (the active/passive marks space) ``similar'' if all the marks in the two realizations occur at same time and their attribute, active/passive, is the same unless a mark occurs at a time $t$ at a pair $(x,x+1)$ such that $x=x_i(t^-)$, $x+1=x_j(t^-)$ and both $i,j \le \ell$: in such a case the mark attributes in $\om$ and $\om'$ may either be equal or opposite.  It is readily seen that this is an equivalence relation. Common to all $\om$ in a same equivalence class  are all the times $s_1<\dots<s_{M'}$ in $[0,t]$ where  a mark appears involving sites with at least one stirring particle with label $\le \ell$.  We call $t_1<\dots<t_{M}$ its subset when both sites indicated by the mark are occupied by particles with label $\le \ell$. We define $\delta(t_i)=\pm 1$ if the mark at time $t_i$ is active, respectively passive. Then, conditioned on the  equivalence class, the variables $\delta(t_i)$ are independent with probability $\frac 12, \frac 12$.

  We define
  $$D(s)=\big(d_1(s),\dots,d_\ell(s)\big),\quad  d_i(s):= x_i(s)-x^0_i(s),\quad s\le t$$
and we observe that for a mark $\{x,x+1\}$ at time $s_i\notin \{t_1,.,t_{M}\}$ we have  $D(s_i^+)=D(s_i^-)$ if $\{x,x+1\}\cap\{-N,N\}=\emptyset$; if instead $\{x,x+1\}\cap\{-N,N\}\neq\emptyset$, then $|d_k(s_i^+)|\le |d_k(s_i^-)|$ for all $k$.

We call $L(s_i)$ the label of the particle involved by the mark at time $s_i$ if $s_i\notin \{t_1,.,t_{M}\}$, otherwise
$L(s_i)$ is the largest of the two labels involved. $L(s_i)$ is specified by the values of all the $\delta(t_k)$ with $t_k<s_i$.
Let $\alpha(t_i)\in \{-1,1\}$ (which depends on all the previous history) be
such that $d_{L(t_i)}(t_i^+)-d_{L(t_i)}(t_i^-)=\alpha(t_i)\delta(t_i)$, namely  if $\{x,x+1\}$ is the mark at time $t_i$, then $\alpha(t_i)=1$ ($=-1$) if $x_{L(t_i)}(t_i^-)=x$ ($=x+1$ resp.).

We now define an auxiliary process $D^*=\big(d^*_1(s),\dots,d^*_\ell(s)\big)$, $s\le t$,  that jumps only at the times $t_i$ when only the $L(t_i)$ component varies (by $\pm 1$). Suppose we have specified the values $\delta(t_j)$, $j<i$, so that we know $D(s)$, $s<t_i$, and suppose inductively that we know $D^*(s)$ as well.
Let $e_k$ denote the vector $\mathbf 1_{i=k}, i=1,\dots,\ell$.  If $k=L(t_i)$ we define
                      \begin{equation}
              \label{19.e12}
              d_k^*(t_i^+)-d_k^*(t_i^-)=  \beta(t_i)\delta(t_i), \quad
D^*(t_i^+)-D^*(t_i^-)=  e_{L(t_i)} \beta(t_i)\delta(t_i)
            \end{equation}
where $\beta(t_i)= \alpha(t_i)$ if  $d_k^*(t_i^-)$ and $d_k(t_i^-)$ have same sign otherwise  $\beta(t_i)= -\alpha(t_i)$ (with the convention that two numbers have same sign also when one of them is 0).
The above defines inductively $D^*(s)$  at all times having supposed $D^*(0)=0$.

From the definition we can inductively check that $|d_i(s)| \le |d^*_i(s)|+1$ for all $s\le t$ and all $i\in \{1,\dots,\ell\}$. Indeed it is enough to prove that if this holds up to  $t_j^-$ then it is
also true at $t_j^+$. If
$d_i(t_j^-)$ and $d^*_i(t_j^-)$ have the same sign (including the case when one or both are 0) then  $d_i(t_j^+)-d_i(t_j^-)=d^*_i(t_j^+)-d^*_i(t_j^-)$ and so the inequality holds by the induction hypothesis. If instead they have opposite sign  $d_i(t_j^+)-d_i(t_j^-)=-\big(d^*_i(t_j^+)-d^*_i(t_j^-)\big)$ and again the inequality holds by the induction hypothesis.
We have
          \begin{equation}
              \label{19.e13}
D(t)^2 \le 2D^*(t)^2 +2\ell
            \end{equation}
By \eqref{19.e12}
         \begin{equation}
              \label{19.e14}
D^*(t_M^+)^2 = D^*(t_M^-)^2 + 2\beta(t_M)\delta(t_M) \langle D^*(t_M^-), e_{L(t_M)} \rangle + 1.
            \end{equation}
Denoting by $\mathcal P$ and $\mathcal E$ law and  expectation  conditioned to an equivalence class and denoting by $\mathcal F_{t_M^-}$ the $\si$-algebra generated by all $\delta(t_i)$ with $i<M$ we have using \eqref{19.e14}
                      \begin{equation}
              \label{19.e15}
\mathcal E\Big[|D^*(t_M^+)|^{2n} \big| \mathcal F_{t_M^-}\Big] = |D^*(t_M^-)|^{2n}+\sum_{m_1+m_2\le n, m_2\ne 0\; {\rm even}} c_{m_1,m_2} |D^*(t_M^-)|^{2m_1}  (\langle D^*(t_M^-), e_{L(t_M)} \rangle)^{m_2}
             \end{equation}
the sum being over even $m_2$ because $\beta(t_M)$ and $ D^*(t_M^-)$ are $\mathcal F_{t_M^-}$-measurable and
   $$ \mathcal E [\delta(t_M) | \mathcal F_{t_M^-}] =\mathcal E [\delta(t_M)]= 0$$
because the $\delta(t_i)$ are independent.
Bounding $|\langle D^*(t_M^-), e_{L(t_M)} \rangle|^{m_2} \le |D^*(t_M^-)|^{m_2}$ and recalling that $m_2$ is even,  we get (for suitable coefficients $c_k$)
                      \begin{equation}
              \label{19.e16}
\mathcal E\Big[\,|D^*(t_M^+)|^{2n} \big| \mathcal F_{t_M^-}\Big] \le |D^*(t_M^-)|^{2n}+ \sum_{k=1}^{n-1} c_k| D^*(t_M^-)|^{2k}
             \end{equation}
and, by iteration,
                    \begin{equation}
              \label{19.e17}
\mathcal E\Big[ |D^*(t)|^{2n}\Big] \le \sum_{m=1}^{M}\sum_{k=1}^{n-1} c_k\mathcal E\Big[ | D^*(t_m^-)|^{2k}\Big].
             \end{equation}
Then, by induction in $n$,  there are new coefficients $c_n$ so that
$$
\mathcal E\Big[ |D^*(t)|^{2n}\Big] \le c_n M^n
$$
\eqref{19.e11} then follows using the Chebyshev inequality and \eqref{19.e2}, details are omitted. \qed

\vskip2cm

\section{Integral inequalities for the $v$-functions}
    \label{sec:3}
\vskip.5cm

                  \noindent
The difference between the true expectation of $\eta(x,t)$ and $\rho_\eps(x,t)$ will be controlled by the  $v$-functions:


\vskip.5cm

{\bf Definition.} [{\em $v$ functions}]   We fix arbitrarily $\eta \in \{0,1\}^{\La_N}$, shorthand by $\E_\eps$ the expectation for the  process with generator $L_\eps=\eps^{-2}L_0 +\eps^{-1} L_b$ which starts from  $\eta$ and write $\rho_\eps(x,t)$ for the solution of \eqref{pro.2.1}
with initial condition $\eta$. Recall  that the $v$-functions are defined in \eqref{pro.2.3.1} and, for brevity, we shall write $v(\und x,t)\equiv v^\eps(\und x,t|\eta)$.

              \vskip1cm

\noindent
{\bf Definition.}\; [{\em The $A$, $B$ and $C_\eps$ operators}] These are linear operators acting on $v$ functions.
For any $X \subset \La_N$,  $t>0$ we define $(Av)(X,t)=0$ if $|X|<2$ while for $|X|\ge 2$ we set:
   \begin{eqnarray}
(Av)(X,t) & := & \sum_{x,y\in X}\mathbf 1_{|x-y|=1}\left\{
\[
\rho_\eps(x,t)-\rho_\eps(y,t)
\]
\[
v(X\setminus x,t)-v(X\setminus y,t)
\]\right .\nonumber\\
&& \left. -\frac 12 \[ \rho_\eps(x,t)-\rho_\eps(y,t)
\]^2 v(X\setminus(x\cup y),t)
\right\}.\label{11}
       \end{eqnarray}
Given real numbers $b(Z,Z',t)$  with $Z$ and $Z'$ either both subsets of $I_+$ or both subsets of $I_-$,
we define
       \begin{equation}
       \label{12}
(B_{\pm}v)(X,t):=\sum_{Z'\subset I_{\pm}}
b\Big([X\cap I_{\pm}],Z',t\Big) v\Big(X\setminus [X\cap I_{\pm}]\cup Z',t\Big)
       \end{equation}
and set $Bv:=B_+v+B_-v$.  We shall shorthand
       \begin{equation}
       \label{12aa}
 (C_{\eps}v)(X,t):=
\eps^{-2}(Av)(X,t)+\eps^{-1}(Bv)(X,t).
       \end{equation}

              \vskip1cm

               \begin{lemma}
               \label{l1}
For any $X\subset\La_N$ and any $t\geq 0$,
       \begin{equation}
       \label{10}
\frac{d}{dt} v(X,t)=\eps^{-2} (L_0 v)(X,t)+(C_{\eps}v)(X,t),\,\,\,\,\,
     \end{equation}
where
$L_0v$ denotes the action of $L_0$  on  $ v(\cdot,t)$  as a  function of  $X$  with $X$ regarded as a particle configuration. $(C_{\eps}v)$ is defined in \eqref{12aa} with
coefficients $b(Z,Z',t)$  having the following property:
       \begin{eqnarray}
       \label{12ab}
 && b(\emptyset, Z',t)=0,\;\;\; b(Z,\emptyset,t)=0\,\; {\rm if }\, |Z|=1
 \\&& \text{for any integer $M$:}\,\, \sup_{t,|Z|\le M, |Z'|\le M}|b(Z,Z',t)| <\infty.\nn
       \end{eqnarray}

    \end{lemma}

    \vskip.5cm
\noindent
{\it Proof.}
We obviously have
              \begin{equation*}
\frac{d}{dt} v(X,t)=\E_\eps\[\;L_\eps
\prod_{x\in X}\{
\eta(x,t)-\rho_\eps(x,t)\} \]+ \E_\eps\[\;\frac{\partial}{\partial t}
\prod_{x\in X}\{
\eta(x,t)-\rho_\eps(x,t)\}
\]
              \end{equation*}
where the partial derivative acts only on $\rho_\eps(\cdot,t)$. Recalling \eqref{pro.2.1}, when the time derivative acts on the factor $\rho_\eps(x,t)$ it gives rise to the sum $\eps^{-2}\Delta \rho_\eps(x,t)+ \eps^{-1}D\rho_\eps(x,t)$ with $D\rho_\eps=\mathbf 1_{x\in I_+} D_+\rho_\eps(x,t)$ $
-\mathbf 1_{x\in I_-}  D_-\rho_\eps(x,t)$. All terms with $\eps^{-2}\Delta \rho_\eps$ combined with those arising from  $\eps^{-2}L_0$ are the same as when  $L_b$ is absent and it is proved
in Lemma 10.1.2 of \cite{DP} that their sum is equal to $\eps^{-2}[L_0 v + Av]$. Thus we need only to prove that the remaining terms (arising from the action of $L_b$ and from the terms with $\eps^{-1}D\rho_\eps$)
is given by \eqref{12}. Considering the  terms arising from the boundary generator in $I_+$ (the one in  $I_-$ is similar and the analysis omitted) we get, modulo a pre-factor  $\eps^{-1}$,
        \begin{equation}
      \label{12.1a}
\E_\eps\[
L_{b,+}
\prod_{x\in X}[
\eta(x,t)-\rho_\eps(x,t)
]
\]
-\sum_{x\in X}
\E_\eps\[
D_+\rho_\eps(x,t)
\prod_{y\in X\setminus x}[
\eta(y,t)-\rho_\eps(y,t)
]
\].
 \end{equation}
We write
              \begin{equation*}
L_{b,+}\eta(x,t) = \frac{j}{2}\Big([1-\rho_\eps(x,t)]-[\eta(x,t)-\rho_\eps(x,t)]\Big)
\prod_{x+1\le y \le  N}\Big([
\eta(y,t)-\rho_\eps(y,t)]+\rho_\eps(y,t)
\Big)
       \end{equation*}
After expanding the products and denoting below by $Z'$ a subset of $\{x+1,\ldots,N\}$,
         \begin{eqnarray*}
&&\hskip-1.6cm L_{b,+}\eta(x,t) =  D_+\rho_\eps(x,t) + \frac{j}{2}(1-\rho_\eps(x,t)) \sum_{Z' \ne \emptyset} \{\prod_{z\notin Z'}\rho_\eps(z,t)\}\{\prod_{y\in Z'}[\eta(y,t)-\rho_\eps(y,t)]\}
 \nonumber\\
&& -\frac{j}{2}[\eta(x,t)-\rho_\eps(x,t)]
\sum_{Z'}
\{\prod_{z\notin Z'}\rho_\eps(z,t)\}\{\prod_{y\in Z'}[\eta(y,t)-\rho_\eps(y,t)]\}
           \end{eqnarray*}
(in the first sum $Z'=\emptyset$ is absent because it has been included in  $D_+\rho_\eps(x,t)$).

\noindent
The term with $D_+\rho_\eps(x,t)$ cancels with the second term in \eqref{12.1a} if
$|X\cap I_+|=1$. 
Hence all remaining terms have at least
one
factor $\eta-\rho_\eps$ in $I_+$. The other properties of  the coefficients $b$
stated in \eqref{12ab} easily follow.
\qed

\vskip.5cm

For the stirring process defined in Section \ref{secN6}, we let $\PP_\eps(X\overset{s}\rightarrow Y)=\PP_\eps(X(s)=Y|X(0)=X)$, for $X\subset \La_N$ to $Y\subset \La_N$, $|X|=|Y|$. (In particular  $\PP_\eps(X\overset{s}\rightarrow Y)=P^{(\eps)}_s(x,y)$ as in \eqref{a3.6}, if $X=\{x\}$, $Y=\{y\}$.) Since $v(X,0)=0$, the integral form
of \eqref{10} is
         \begin{equation}
         \label{14}
v(X,t)=\int_0^t\,ds
\sum_Y\PP_\eps(X\overset{s}\rightarrow Y)(C_{\eps}v)(Y,t-s).
         \end{equation}
We start by bounding the contribution of $B_{\pm}v$ to $C_\eps v$ in the right hand side of \eqref{14}:

\vskip.5cm

                 \begin{lemma}
                 \label{l2}
For any $n$ and any $\zeta>0$ there is a constant $c$ so that for any $X\subset\La_N$, $|X|=n$, and any $s<t \le \log \eps^{-1}$

              \begin{eqnarray}
              \label{13}
&& \sum_Y\PP_{\eps}(X\overset{s}\rightarrow Y) \eps^{-1}|(B_{\pm}v)(Y,t-s)|
\leq
\sum_{Z'\subset I_{\pm}}
\sum_{\emptyset \ne X''\subset X} \big(1-\mathbf 1_{|Z'|=0,|X''|=1}\big) \nn\\&&\hskip2cm \times
\frac{c\eps^{-1}}{(\eps^{-2}s)^{|X''|/2}+1}
 \sum_{W\subset I_{\pm}^c}\PP_\eps(X\setminus X''\overset{s}\rightarrow W) |v(W\cup Z',t-s)|.
  \end{eqnarray}

   \end{lemma}

   \vskip.5cm

\noindent
{\it Proof.}  We consider explicitly only the case with $B_+$. The left hand side of \eqref{13} without absolute values is equal to:
             \begin{equation}
             \label{14.1}
\Psi:=\eps^{-1}\sum_{Y,Z'\subset I_+}\PP_\eps(X\overset{s}\rightarrow Y)\mathbf 1_{Y\cap I_+\neq\emptyset}\,
b(Y\cap I_+,Z') v\Big([Y\setminus (Y\cap I_+)]\cup Z',t-s\Big).
              \end{equation}
We decompose $Y=W\cup Z$, $W\subset I_+^c$ and $Z \subset I_+$, $|Z|>0$. Then  the sum over $Y$ becomes a sum over $W$ and $Z$ with the condition that $|W\cup Z|=|X|$. For each fixed $W$ and $Z$ we apply
Andjel's inequality, see \cite{andjel}, and get
        \begin{equation*}
\PP_\eps(X\overset{s}\rightarrow W\cup Z)\leq   \{\sum_{Y\supset W}\PP_\eps(X\overset{s}\rightarrow Y)\}\;\{\sum_{Y'\supset Z}\PP_\eps(X\overset{s}\rightarrow Y')\}.
         \end{equation*}
By Liggett's inequality, see \cite{liggett},
              \begin{equation*}
\PP_\eps(X\overset{s}\rightarrow Y')\leq  \prod_{x\in X} \sum_{y\in Y'}\PP_\eps(x\overset{s}\rightarrow y)
              \end{equation*}
and by \eqref{N4.5}
there is $c=c_{|X|}$ so that
              \begin{equation*}
\sum_{Y'\supset Z}\PP_\eps(X\overset{s}\rightarrow Y')\leq  \frac{c}{(\eps^{-2}s)^{|Z|/2}+1}.
              \end{equation*}
Observe that
 $$\sum_{Y\supset W}\PP_\eps(X\overset{s}\rightarrow Y)=\sum_{X'\subset X}\PP_\eps(X'\overset{s}\rightarrow W),$$
and collecting the estimates we have
  \begin{eqnarray*}
&&\hskip-1cm|\Psi|\le  \sum_{Z \subset I_+, |Z|\ge 1} \sum_{W\subset I_+^c, |W\cup Z|=|X|}\frac{c\, \eps^{-1}}{(\eps^{-2}s)^{(|Z|)/2}+1}
\sum_{Z'\subset I_+}|b(Z,Z') |
 \\&&\hskip2cm \times \sum_{X'\subset X, |X'|=|W|}
\PP_\eps(X'\overset{s}\rightarrow W)|v\big(W\cup Z',t-s\big)|.
       \end{eqnarray*}
From  \eqref{12ab} we get that $|b(Z,Z') |\le c \big( 1-\mathbf 1_{|Z|=1,|Z'|=0}\big)$.
Thus, denoting by $p(K)$ the number of subsets of $I_+$,
        \begin{eqnarray*}
&&\hskip-1cm|\Psi|\le  \sum_{\emptyset \ne X''\subset X}
\sum_{Z'\subset I_+}(1-\mathbf 1_{|X''|=1,|Z'|=0})
 \frac{c\, p(K)\eps^{-1}}{(\eps^{-2}s)^{(|X''|)/2}+1} \\&&\hskip2cm \times\sum_{W\subset I_+^c:|W|= |X\setminus X''|}
\PP_\eps(X\setminus X''\overset{s}\rightarrow W)| v\big(W\cup Z',t-s\big)|.
       \end{eqnarray*}
 \qed

\vskip.5cm

We are left with the bound  of the contribution in  \eqref{14} due to $\eps^{-2}A$, see \eqref{11}.
It is crucial here  to exploit the smallness of the gradients, namely the differences $\rho_\eps(x,t)-\rho_\eps(y,t)$ and
$v(X\setminus x,t)-v(X\setminus y,t)$ (recall that $|x-y|=1$). Both bounds use the parabolic nature of the evolution, but the latter
requires a more delicate analysis which, following \cite{DP}, is based on the realization of the stirring process given in Definition 1 of Section \ref{secN6}.

Recalling Definition 2 in Section \ref{secN6}, we order arbitrarily the sites of $X$ 
 which are then denoted by $\und x=(x_1,\dots,x_n)$ and set $v(\und x,t) := v(X,t)$,
$v(\und x,t)$ being symmetric under exchange of labels. We denote by $\EE_{\eps,\und x}$ the expectation with respect to the stirring process defined in Section \ref{secN6} starting at time 0 from $\und x$ and when the starting point will be clear from the context we shorthand $\EE_\eps\equiv \EE_{\eps,\und x}$.
We then rewrite \eqref{13} as
             \begin{eqnarray}
              \label{13.1.1}
&& \sum_Y\PP_\eps(X\overset{s}\rightarrow Y) \eps^{-1}|(B_{\pm}v)(Y,t-s)|
\leq \sum_{\emptyset \ne J \subset \{1,\dots,n\}}
\sum_{\und z'\subset I_{\pm}} [1-\mathbf 1_{|J|=1,|\und z'|=0}] \nn\\&&\hskip2cm \times
\frac{c\eps^{-1}}{(\eps^{-2}s)^{|J|/2}+1} \EE_{\eps,\und x}\Big[ \mathbf 1_{\und x^{(J)}(s)\subset   I_{\pm}^c}  |v\big(\und x^{(J)}(s)\cup \und z',t-s\big)|\Big]
  \end{eqnarray}
where $\und x^{(J)}$ is the configuration obtained by erasing from $\und x$ all $x_j$ with $j\in J$ and
we then say that all the particles $x_j, j\in J$ have died (and their labels will not be used again) and that the particles $\und z'$ are born  at time $s$. Our general rule  to label a new particle is to use the smallest integer never used earlier in the labeling (the order in which the particles in $\und z'$ are born is chosen arbitrarily).  In conclusion \eqref{13.1.1} describes a labeled stirring evolution with a death/birth process at time $s$ (notice that it is the same to erase the particles $x_j, j\in J$ either at time 0 or at time $s$).
In an analogous way we  write the labeled version of \eqref{11} as
   \begin{eqnarray}
(Av)(\und x,t) & := & \sum_{x_i,x_j\in \und x}\mathbf 1_{|x_i-x_j|=1}\Big\{
[
\rho_\eps(x_i,t)-\rho_\eps(x_j,t)
]
[
v(\und x^{(i)},t)-v(\und x^{(j)},t)
] \nonumber\\
&&  -\frac 12  [ \rho_\eps(x_i,t)-\rho_\eps(x_j,t)
 ]^2 v(\und x^{(i,j)},t)
\Big\}
       \label{11.1}
       \end{eqnarray}
so that we have:
   \begin{eqnarray}
&& \hskip-2cm |v(\und x,t)| \le \int_0^t ds \bigg(\sum_{u=\pm}\,\,
\sum_{\emptyset \ne J \subset \{1,\dots,n\}}
\sum_{\und z'\subset I_u}
\frac{c\eps^{-1}}{(\eps^{-2}s)^{|J|/2}+1}[1-\mathbf 1_{|J|=1,|\und z'|=0}] \nn\\&&\hskip4cm \times \EE_\eps\Big[ \mathbf 1_{\und x^{(J)}(s)\subset   I_{u}^c}  |v(\und x^{(J)}(s)\cup \und z',t-s)|\Big]\nn\\&& \hskip2cm + \; \eps^{-2} \EE_\eps\Big[
|(Av)(\und x(s),t-s)|\Big] \bigg),
       \label{11.1n}
       \end{eqnarray}
with $Av$ as in \eqref{11.1}.

The difference of the $\rho_\eps$'s in \eqref{11.1} is bounded by using \eqref{N5.2a}.
Indeed by using a weaker form as in
\eqref{18a} we bound the last term in \eqref{11.1n} by
              \begin{eqnarray}
             \label{13n}
&& \hskip-1cm  \EE_\eps\Big[| \eps^{-2} (Av)(\und x(s),t-s)|\Big] \leq    \sum_{i,j}\EE_\eps\Big[ \mathbf 1_{|x_i(s)-x_j(s)|=1}
 \nn\\&& \times c\eps^{-2}\Big(
\frac{|v(\und x^{(i)}(s),t-s)-v(\und x^{(j)}(s),t-s)|}{(\eps^{-2}(t-s))^{1/2-\zeta}+1}
+ \frac{|v(\und x^{(i,j)}(s),t-s)|}{[\eps^{-2}(t-s)]^{1-\zeta}+1}\Big) \Big].
              \end{eqnarray}

To bound  the $v$-gradients appearing in \eqref{13n} we shall use the following lemma  where we take advantage for the first time of the features of  the  active/passive marks process (analogous estimates are  given in Sect.10.1 of \cite{DP}) .

\vskip.5cm

\noindent
{\bf Definition.} [{\em The stopping time $\tau_{i,j,t_0}$}]

Let $\{\und x(t)\}_{t\ge 0}$ be the labeled process realized in the active/passive marks process, let $i$ and $j$ be the labels of two of its particles and $t_0\ge 0$.
We then define $\tau_{i,j,t_0}$ as the first time $\tau > t_0$
when  {\em(i)} $|x_i(\tau)-x_j(\tau)|=1$  and {\em(ii)} at $\tau$ there is a mark (either active or  passive) between  $x_i(\tau)$ and $x_j(\tau)$; otherwise we set $\tau=\infty$.
When $t_0=0$ we just write $\tau_{i,j}$.

\vskip.5cm

         \begin{lemma}
We have
                     \begin{eqnarray}
                     \label{17.3}
|v(\und x^{(i)},t)-v(\und x^{(j)},t)| &\leq &
\int_0^t ds\; \EE_\eps \Big[\mathbf 1_{\tau_{i,j}\ge \frac s2}\sum_{\und y} \{\PP_\eps(\und x^{(j)}(\frac{s}{2})\overset{s/2}\rightarrow \und y)
\nn\\ &+& \PP_\eps(\und x^{(i)}(\frac{s}{2})\overset{s/2}\rightarrow \und y)\} |(C_\eps v)(\und y,t-s)| \Big].
                   \end{eqnarray}

            \end{lemma}

\vskip.5cm

\noindent
{\em Proof.} Denoting by $\und x(s)$ the process starting from $\und x$ with both particles $i$ and $j$,
                    \begin{eqnarray*}
v(\und x^{(i)},t)-v(\und x^{(j)},t)= \int_0^t ds\; \EE_\eps \Big[ (C_\eps v)\big(\und x(s)\setminus x_i(s),t-s\big)- (C_\eps v)\big(\und x(s)\setminus x_j(s),t-s\big) \Big]
                   \end{eqnarray*}
which is the expectation at time $s$ of the function $f(\und y)= (C_\eps v)(\und y^{(i)},t-s)-
(C_\eps v)(\und y^{(j)},t-s)$. Since $f$ is antisymmetric under the exchange of particles $i$ and $j$,
\eqref{17.3} follows from Lemma \ref{x}.  \qed

\vskip.5cm
The reason for the time interval $s/2$ in the lemma is to be able to exploit   \eqref{13.1.1}. We have in fact from \eqref{17.3} and \eqref{13.1.1}
                     \begin{eqnarray}
                     \label{17.3.3}
&&\hskip-2cm  |v(\und x^{(i)},t)-v(\und x^{(j)},t)| \leq
\int_0^t ds\; \bigg( \sum_{u=\pm}\;\;\sum_{\emptyset \ne J \subset \{1,\dots,n\}}\frac{c\eps^{-1}}{(\eps^{-2}\frac s2)^{|J|/2}+1}
\sum_{\und z'\subset I_{u}}[1-\mathbf 1_{|J|=1,|\und z'|=0}]
 \nn
\\&&\hskip3cm   \times\EE_\eps \Big[\mathbf 1_{\tau_{i,j}>\frac s2} \mathbf 1_{\und x^{(J)}(s)\in I_u^c} |v(\und x^{(J)}(s)\cup \und z',t-s)| \Big]
\nn \\ &+& \EE_\eps\Big[ \mathbf 1_{\tau_{i,j}>\frac s2}\eps^{-2}\{|(A v)(\und x^{(i)}(s),t-s)|+|(A v)(\und x^{(j)}(s),t-s)|\} \Big]\bigg),
                   \end{eqnarray}
while the $v$ functions alone are bounded by:
             \begin{eqnarray}
              \label{13.1.111}
&& |v(\und x,t)|\le c\int_0^t ds \Big( \sum_{u=\pm}\;\; \sum_{J\subset \{1,\dots,n\}} \frac{\eps^{-1}}{(\eps^{-2}s)^{|J|/2}+1}\sum_{ \und z' \subset I_u} [1-\mathbf 1_{|J|=1,|\und z'|=0}]
 \nn
\\&&\hskip4cm   \times
\EE_\eps \Big[ \mathbf 1_{\und x^{(J)}(s)\subset I_u^c} |v\big( \und x^{(J)}(s)\cup \und z',t-s\big)|\Big] \nn\\&& \hskip2cm  +\sum_{i,j\subset \{1,\dots,n\}}\EE_\eps\Big[ \mathbf 1_{|x_i(s)-x_j(s)|=1}\,\, \eps^{-2}\Big(
\frac{|v(\und x^{(i)}(s),t-s)-v(\und x^{(j)}(s),t-s)|}{(\eps^{-2}(t-s))^{1/2-\zeta}+1}
\nonumber\\
&& \hskip6cm  +
\frac{|v(\und x^{(i,j)}(s),t-s)|}{[\eps^{-2}(t-s)]^{1-\zeta}+1}\Big).
\Big]
  \end{eqnarray}
We shall derive the desired bound for $|v(\und x,t)|$ by   iterating \eqref{13.1.111} and using \eqref{17.3.3} (complemented by \eqref{11.1} to write the terms $Av$ whenever a  $v$-gradient appears. The series obtained in this way is described in the next section and studied in Section \ref{s5} and  Section \ref{s5a}.

\vskip2cm

\section{The truncated hierarchy}
    \label{sec.4}

By iterating \eqref{17.3.3}--\eqref{13.1.111} we can write the solution as a formal series, but we do not know whether it converges.
We shall therefore truncate the expansion proving that at least for small times the remainder is small. Thus, in a first step we only
prove short time estimates:

\vskip.5cm

                \begin{thm}
                \label{t1}
For any $ c^*<\frac{1}{4(K+2)}$ the following holds.
For any $\beta^*>0$, any initial configuration $\eta_0$
and any positive integer $n$ there is $c$ so that
             \begin{equation}
             \label{0}
|v(\und x,t|\eta_0)|\leq c (\eps^{-2}t)^{-c^* n},\,\,\,\,\,\, t\leq\eps^{\beta^*},\;\; |\und x|=n.
             \end{equation}

             \end{thm}

\vskip.5cm
For brevity in the sequel we shall simply write $v(\und x,t)$ for $v(\und x,t|\eta_0)$.
The theorem will be proved in Section \ref{s6}
using the results of Sections~\ref{s5} and \ref{s5a}. Of course we only need to prove \eqref{0} when $t \ge\eps^2$ because for all values of its arguments
$|v(\und x,t)|\leq 1$.

\vskip.5cm

\noindent
{\bf The setup.}
As mentioned above, Theorem \ref{t1} will be proved by finitely many iterations of
the integral inequalities \eqref{17.3.3}--\eqref{13.1.111}.  The  series obtained in this way will be referred to as ``the truncated hierarchy''.  The number of iterations will depend on $n$ and $\beta^*$, and will be denoted by $M$; its actual value will be specified later in \eqref{8e.2.1} and \eqref{20.107.03}. At each iteration the number of ``particles'', i.e.\  elements in the argument of the $v$-function, increases at most by $K-1$ so that the total number of particles is not larger than $n+(K-1)M$. Hence all constants that appear in the previous section which depend on the cardinality of the configuration in the $v$ functions are bounded by a constant (once $n$ and $\beta^*$ are fixed).  The various terms which appear in the expansion will be classified in terms of sequences called skeletons.  We shall first define the skeletons and then establish the correspondence with the terms appearing in the expansion.  The positions of the particles are not recorded in the skeleton; it says which particles are alive at each step of the process, as well as those which die and are born, specifying also the positions of the new-born particles at their birth.

\vskip.5cm

\noindent
{\bf Definition.} [{\em The skeleton}]

Skeletons are denoted by $\pi$. Each $\pi$ consists of a sequence
$\pi=(\pi_i)_{i=1,\dots,m(\pi)}$, $m(\pi)\le M$ (see the ``setup'' paragraph). ``$i$'' is  a``branching time'' and $\pi_i$ describes the nature of the branching (which particles die and which are born). As we shall explain the alive particles at step $i$, denoted by $A_i$, are determined by the values  $\pi_j$ with $j\le i$ while  the ``initially alive particles'' are  $A_0=\{1,\dots,n\}$.

%
%


$\bullet$\; for each $i$, $\pi_i=(\delta_i,J_i,u_i,\und z_i)$, $\delta_i\in\{0,1,2\}$,
$J_i$ is a finite increasing sequence of distinct positive integers, $u_i\in\{0,+,-\}$,  $\und z_i$ is a labeled configuration, its labels will be denoted by $J_i^+$.
There are several constraints relating the elements $\pi_i$ of  $\pi$ which we state inductively.  We suppose
that we have already chosen the elements $\pi_j$ with $j<i$ and thus know the sequence $A_j$, $j<i$, of alive particles at the branching times $j$. We then want to specify the possible values of $\pi_i$, and for each  choice of $\pi_i$ we shall define $A_i$. (When $i=1$ we only need $A_0$ which is $\{1,\dots,n\}$, hence the induction is complete).

With a small abuse of notation we sometimes identify  $J_i$ with the set of its constituents.

$\bullet$\; Before entering into all the details we just say that the particles which die are those with labels in $J_i$ except when $\delta_i=0$ in which case $J_i$ consists of two particles but only one of them dies; in all cases
$J_i\subset A_{i-1}$. If $u_i=0$ then no particle is born and $\und z_i=\emptyset$.
If $u_i\ne 0$ there may be new particles. The configuration of the new particles is  $\und z_i$ which is contained in $I_{u_i}$. The  labels in   $\und z_i$ are  $J_i^+$, $J_i^+$ is a sequence of consecutive integers, the first one is $h+1$ if $h$ is the max over all integers in  the union of $A_j$ over $j<i$.
The positions of the particles in $\und z_i$ are increasing functions of the labels.

$\bullet$\; If $\delta_i=0,1$ then  $J_i$ is an ordered pair, $J_i=(k_i,\ell_i)$,  $k_i<\ell_i$
$u_i=0$ and $\und z_i=\emptyset$.  If $\delta_i=1$,  then both particles die, so that $A_i=A_{i-1}\setminus \{k_i,\ell_i\}$.  If $\delta_i=0$, then only particle  $\ell_i$ dies so that $A_i=A_{i-1}\setminus \ell_i$.

$\bullet$\; If $\delta_i=2$,   $J_i$ is  non empty ($|J_i|\ge 1$), $u_i\ne 0$ and  $A_i=(A_{i-1}\setminus J_i)\cup J_i^+$
if $J_i^+$ are the labels of the particles in $\und z_i$.  The configuration $\und z_i$ may also be empty but not when
$|J_i|= 1$, in which case $|J_i^+|>0$. 

$\bullet$\; Finally   if $m(\pi)=m<M$ then $\delta_m>0$,
$\und z_m=\emptyset$ and $|J_m|= |A_{m-1}|$, i.e.\ $A_m=\emptyset$. If on the other hand $m(\pi)=M$ then there is no restriction on $A_M$, which could be either empty or non empty.

\vskip.5cm

\noindent
{\bf Definition.} [{\em The branching process}]

Given an element $\om\in \Om$,  the active/passive marks space, $\und x$ (the initial configuration), a skeleton $\pi$ and
$0=t_0 <t_1< \dots <t_m<t_{m+1}=t$, $m=m(\pi)$, we define  $\und x(t)$ by following in the time intervals
$(t_i,t_{i+1})$ the active/passive marks. At time $t_i$ all particles $x_j(t_i^-)$, $j\in J_i$, disappear from $\und x(t_i^-)$ except when
$\delta_i=0$: in that case the particle with label $k_i$ remains alive, the one with label $l_i$ survives but it will die at time $t_i+(t_{i+1}-t_i)/2$. We also require that if $\delta_i=2$ then $x_j(t_i^-)\in I_{u_i}^c$ for  all $j \in A_{i-1}\setminus J_i$, we shall shorthand this event by $R_i$.
We complete the definition of $\und x(t)$ by saying that at time $t_i^+$ we add the labeled particles  $\und z_i$.

\vskip.5cm

\noindent
In order to write ``the truncated hierarchy'' we introduce the factors  $\ga_i$ which depend on the realization of the active/passive marks process, $\und x$ (the initial configuration), the skeleton $\pi$ and the sequence of times
$0=t_0 <t_1< \dots <t_m<t_{m+1}=t$, $m=m(\pi)$.

\begin{itemize}

\item If $\delta_i=0$ then
              \begin{equation}
              \label{19.1}
   \ga_i= \mathbf 1_{|  x_{k_i}(t_i)- x_{l_i}(t_i)|=1} \frac{\eps^{-2}}{[\eps^{-2}(t-t_i)]^{\frac 12-\zeta}+1}
   \mathbf 1_{\tau_{k_i,l_i,t_i}>t_i+\frac{t_{i+1}-t_i}{2}}
             \end{equation}
   It means that we are considering the first term in the second expectation on the right hand side of \eqref{13.1.111} with $i,j$ equal to $(k_i,l_i)$ and then, when writing the $v$-gradient via \eqref{17.3.3}, we take the term in the second expectation where the label $l_i$ is missing (i.e.\ particle $k_i$ survives, particle $l_i$ dies).

\item If $\delta_i=1$ then
              \begin{equation}
              \label{19.2}
   \ga_i= \mathbf 1_{|  x_{k_i}(t_i)- x_{l_i}(t_i)|=1} \frac{\eps^{-2}}{[\eps^{-2}(t-t_i)]^{1-\zeta}+1}
             \end{equation}
  It means that we are considering the second term in the second expectation on the right hand side of \eqref{13.1.111} with $(i,j)$ equal to $(k_i,l_i)$.

\item If $\delta_i=2$ then
               \begin{equation}
              \label{19.3}
   \ga_i=  \mathbf 1_{ R_i} \frac{\eps^{-1}}{[\eps^{-2}(t_i-t_{i-1})]^{p_i/2}+1},\quad p_i\equiv |J_i|
             \end{equation}
recalling that $R_i=\{ x_j(t_i^-)\in I^c_{u_i}, j\notin J_i\}$.
In \eqref{19.3} we are considering the first term on the right hand side of \eqref{13.1.111} with $u=u_i$, $J=J_i$ and $\und z'=\und z_i$.

\end{itemize}

\vskip1cm
In this way we have classified all possible terms of the truncated hierarchy and
              \begin{equation}
              \label{19}
|v(\und x,t)|\leq c\sum_{\pi} w_{\pi}(\und x,t)
          \end{equation}
with $c=c(n,\beta)$ is a constant (as discussed in the setup definition) and $w_{\pi}(\und x,t)$ is obtained by integrating the product of all the $\ga_i$ defined above:
             \begin{eqnarray}
             \label{20}
&&\hskip-1.6cm       w_{\pi}(\und x,t) =
\int_0^t dt_1\ldots\int_{t_{m-1}}^t dt_m
\prod_{\delta_i=1}\[
\frac{\eps^{-2}}{[\eps^{-2}(t-t_i)]^{1-\zeta}+1}
\]
\prod_{\delta_i=0}\[
\frac{\eps^{-2}}{[\eps^{-2}(t-t_i)]^{1/2-\zeta}+1}
\]\nonumber\\
&&\hskip.3cm
\prod_{\delta_i=2}\[
\frac{\eps^{-1}}{[\eps^{-2}(t_i-t_{i-1})]^{p_i/2}+1}\]
\EE_\eps \[
\prod_{\delta_i=0,1}
\mathbf 1_{|x_{k_i}(t_i)-x_{l_i}(t_i)|=1}
\prod_{\delta_i=0}\mathbf 1_{T_i}  \prod_{\delta_i=2}  \mathbf 1_{R_i}
\]
            \end{eqnarray}
where $m=m(\pi)$, $p_i$ is defined in \eqref{19.3}, the product over $\{\delta_i=k\}$,  $k=0,1,2$,  means the product over $\{i\in\{1,\ldots,m\}:\,\delta_i=k\}$ and
              \begin{equation}
              \label{19.4}
T_i:=\{\tau_{k_i,l_i,t_i}>t_i+\frac{t_{i+1}-t_i}{2}\},\;\; R_i:=\{ x_j(t_i^-)\in I^c_{u_i}, j\notin J_i\}
             \end{equation}
The expectation $\EE_\eps$ is with respect to the active/passive marks process and $\und x(t)$ is the branching process defined above (in terms of $\pi$ and of the realization of  the active/passive marks process).  If $m(\pi)=M$ there could be surviving particles at $t_M$, i.e.\ a $v$-function $|v(\und x(t_M),t-t_M)|$ that in \eqref {20} has been  bounded by $1$.

\vskip2cm

   \section{Bounds when times do not cluster}
   \label{s5}
The proof of Theorem~\ref{t1} is based on bounds of $w_\pi(\und x,t)$ which will be proved in this and in the next sections. As in the statement of Theorem~\ref{t1} $n$, the cardinality of $\und x$, is fixed and $t\le \eps^{\beta^*}$ $\beta^*>0$.  As already mentioned in Section \ref{sec.4} we then introduce a parameter $M$ which depends on $n$ and
$\beta^*$ and we only  consider  skeletons $\pi$ such that $m(\pi)\le M$. The choice of $M$ will be specified in \eqref{8e.2.1} and \eqref{20.107.03}.  Hereafter $n$, $M$ and $\beta^*$ are to be considered fixed and any parameter which depends only on $n$, $M$ and $\beta^*$ will be called constant.
Setting
    \begin{equation}
               \label{6e.0.0}
a= \frac{K}{K+1}
       \end{equation}
we introduce the quantity $\Delta=\Delta(a,t)$  as
               \begin{equation}
               \label{6e.0.1}
 \Delta = \begin{cases} \eps^a &{\rm if}\; t > (M+1)\eps^a \\\frac{t}{M+1} &{\rm if}\; t \le (M+1)\eps^a .
 \end{cases}
       \end{equation}
The choice of $a$ will be explained later in the course of the proofs.
The parameter $\Delta$ is used to distinguish cases when the times $t_1,\dots,t_m$ ``cluster to $t$'' or not, i.e.\ if $t_m \geq t-\Delta$ or $t_m < t-\Delta$. We accordingly split \eqref{20} writing
              \begin{equation}
               \label{6e.0.2}
 w_\pi(\und x,t)=  w'_\pi(\und x,t) +  w''_\pi(\und x,t)
       \end{equation}
where $ w'_\pi(\und x,t)$ is defined by the right hand side of  \eqref{20} with the integral over $t_m$, $m=m(\pi)$, restricted to $\{t_m<t-\Delta\}$; $w''_\pi(\und x,t)$ is instead the integral over $\{t_m\geq t-\Delta\}$.  The analysis of both
$w'_\pi(\und x,t)$ and  $ w''_\pi(\und x,t)$ consists of two steps: we first apply the results of Section \ref{secN6}
to bound the expectation in \eqref{20}; after this we are reduced to a rather explicit integral over $t_1\dots t_m$ which is bounded in the second  step. Convergence problems in the latter motivate the type of inequalities used in the first step.  The case $\{t_m<t-\Delta\}$ is much simpler and examined first in this section where we shall prove:

\vskip1cm

     \begin{prop}
     \label{prop6.1}
For any $\zeta>0$ there is $c$   so that for all $\pi:m=m(\pi) \le M$, for all
$\und x: |\und x|=n$, for all $\eps>0$ and all
$t \le \eps^{\beta^*}$:
          \begin{equation}
         \label{6e.0.3}
w'_\pi(\und x,t)
\leq c (\eps^2t)^{-\zeta M} \Delta^{-S_1(m)} \eps^{S_2(m)} t^{S_3(m)}
       \end{equation}
 where, recalling that $p_i=|J_i|$ as specified by $\pi$,
            \begin{eqnarray}
         \label{6e.0.4}
&& \hskip-.7cm S_1(m) = |\{i \le m:\delta_i=1\}| + \frac 12 |\{i \le m:\delta_i=0\}| \nn\\&&
\hskip-.7cm S_2(m) = |\{i \le m:\delta_i=0,1\}| +   |\{i\le m:\delta_i=2,\, p_i\geq 2, \delta_{i-1}\neq 0\}| \nn\\&&
\hskip-.7cm S_3(m)=  \frac 12 \Big( |\{i \le m:\delta_i=1\}| +  |\{i \le m:\delta_i=2, p_i=1\}| +   |\{i\le m:\delta_i=2, p_i\geq 2,\delta_{i-1}= 0\}| \Big)\nn\\
       \end{eqnarray}

    \end{prop}

\vskip1cm

The proof to Proposition \ref{prop6.1} is given after stating and proving Lemma \ref{lemma6e.1} below.
We fix $\pi$, write $m=m(\pi)$, $t_0:=0$, $t_{m+1}:=t$ and, with the sets $T_i$ and $R_i$ defined in \eqref{19.4}, we set
            \begin{eqnarray}
    &&\psi_h : =
\prod_{i\le h:\delta_i=0,1}
\mathbf 1_{|x_{k_i}(t_i)-x_{l_i}(t_i)|=1}
\prod_{i\le h:\delta_i=0}\mathbf 1_{T_i}
\prod_{i\le h:\delta_i=2}  \mathbf 1_{R_i},\;\; 1\le h \le m
\nonumber\\
         \label{6e.1}
         \\
&& \phi_h:=
\prod_{i>h:\,\delta_i=0,1}
\frac{1}{[\eps^{-2}(t_i-t_{i-1})]^{1/2}+1}
\prod_{i>h:\,\delta_i=0}
\frac{1}{[\eps^{-2}(t_{i+1}-t_{i})]^{1/2}+1},\;\; 0\le h < m,\nn
\end{eqnarray}
letting $\phi_m:=1$ and $\psi_0=1$.  Then:

\vskip.5cm
             \begin{lemma}
             \label{lemma6e.1}
There is $c$ so that for any $h\le m$
               \begin{equation}
               \label{6e.2}
  \phi_h \EE_\eps \[ \psi_h\] \le c \phi_{h-1} \EE_\eps \[ \psi_{h-1}\].
       \end{equation}

               \end{lemma}

\medskip
{\em Proof.} Let
            \begin{equation}
            \label{6e.3}
\mathcal F(t):= \text{  $\sigma$-algebra generated by the active/passive marks in the interval $[0,t)$  }
              \end{equation}
Suppose first $\delta_h=0$. We then condition on $\mathcal F(t_h)$ getting
               \begin{equation}
               \label{6e.2.0}
 \EE_\eps \[ \psi_h\] = \EE_\eps \[ \psi_{h-1} \mathbf 1_{| x_{k_h}(t^-_h)- x_{l_h}(t^-_h)|=1}\Big(  \PP_\eps\[T_h|\mathcal F(t_h)\]\Big)\]
       \end{equation}
By \eqref{19.e1}
               \begin{equation}
               \label{6e.2.01}
\PP_\eps\[T_h|\mathcal F(t_h)\] \le  \frac{c}{[\eps^{-2}(t_{h+1}-t_h)]^{1/2}+1}.
       \end{equation}
We shorthand $t_{h,+}:= t_h+ \frac{t_{h+1}-t_h}{2}$ and have
               \begin{equation}
              \label{6e.4}
 \EE_\eps \[ \psi_{h-1} \mathbf 1_{| x_{k_h}(t^-_h)- x_{l_h}(t^-_h)|=1} \] = \EE_\eps \[ \psi_{h-1}\Big(  \PP_\eps\[ | x_{k_h}(t^-_h)- x_{l_h}(t^-_h)|=1 | \mathcal F(t_{h-1,+})\]\Big)\].
       \end{equation}
By  \eqref{19.e6}
             \begin{equation}
             \label{6e.5}
 \PP_\eps\[ | x_{k_h}(t^-_h)- x_{l_h}(t^-_h)|=1 \;|\; \mathcal F(t_{h-1,+})\]
\leq\frac{c}{[\eps^{-2}(t_{h}-t_{h-1})]^{1/2}+1}
             \end{equation}
thus completing the proof of \eqref{6e.2} when $\delta_h=0$.  When $\delta_h=1$,
 \eqref{6e.2}  follows from  \eqref{6e.4}--\eqref{6e.5},  while when $\delta_h=2$ we simply bound
 $ \mathbf 1_{R_h} \le 1$.  \qed

\vskip1cm

{\bf Proof of Proposition \ref{prop6.1}.}

Recalling the definitions \eqref{6e.0.2} and \eqref{6e.1},  we apply  repeatedly  Lemma \ref{lemma6e.1} to get
              \begin{equation}
               \label{6e.6}
 w'_\pi(\und x,t)\le  \int_0^t dt_1\ldots\int_{t_{m-1}}^{t-\Delta} dt_m f_{1,\ldots,m}(t_1,\ldots,t_m)
       \end{equation}
where
            \begin{eqnarray}
            \label{20.7}
f_{1,\ldots,m} & := &
\prod_{\delta_i=1}\[
\frac{\eps^{-2}}{[\eps^{-2}(t-t_i)]^{1-\zeta}+1}
\frac{1}{[\eps^{-2}(t_i-t_{i-1})]^{1/2}+1}
\]
\prod_{\delta_i=2}
\frac{\eps^{-1}}{[\eps^{-2}(t_i-t_{i-1})]^{p_i/2}+1}
\nonumber\\
&&
\prod_{\delta_i=0}\[
\frac{\eps^{-2}}{[\eps^{-2}(t-t_i)]^{1/2-\zeta}+1}
\frac{1}{[\eps^{-2}(t_{i+1}-t_{i})]^{1/2}+1}
\frac{1}{[\eps^{-2}(t_i-t_{i-1})]^{1/2}+1}
\]
           \end{eqnarray}

We  bound the factors on the right hand side of \eqref{20.7} as follows. Since $t-t_i>\Delta$ for all $i$ we have
             \begin{equation}
             \label{20.102}
\frac{\eps^{-2}}{[\eps^{-2}(t-t_i)]^{1-\zeta}+1}\leq
\eps^{-2\zeta}\Delta^{-1}
\,\,\,\,\,\text{and}\,\,\,\,\,
\frac{\eps^{-2}}{[\eps^{-2}(t-t_i)]^{1/2-\zeta}+1}\leq
\eps^{-2\zeta}\eps^{-1}\Delta^{-1/2},
              \end{equation}
where we have used that $\Delta^{\zeta}\le 1$  since $\Delta\le \eps^a$, $a>0$.  Moreover we obviously have
        \begin{equation}
        \label{20.1025}
\frac{1}{[\eps^{-2}(t_{i+1}-t_i)]^{1/2}+1}\leq\frac{\eps}{(t_{i+1}-t_i)^{1/2}}.
          \end{equation}
When $\delta_i=2$, $p_i\geq 2$ and $\delta_{i-1}\neq 0$, for any $\zeta>0$ we bound
            \begin{equation}
            \label{20.104}
\frac{\eps^{-1}}{[\eps^{-2}(t_i-t_{i-1})]^{p_i/2}+1}
\leq
\frac{\eps^{-1}}{[\eps^{-2}(t_i-t_{i-1})]^{1-\zeta}}=
\eps^{-2\zeta}\frac{\eps}{(t_i-t_{i-1})^{1-\zeta}}.
               \end{equation}
The first inequality is proved as follows: if $\eps^{-2}(t_i-t_{i-1})\ge 1$   then we drop 1 in the denominator
and replace $p_i/2\ge 1-\zeta$; if $\eps^{-2}(t_i-t_{i-1})\le 1$, then the denominator is
$\ge 1 \ge [\eps^{-2}(t_i-t_{i-1})]^{1-\zeta}$, hence \eqref{20.104}.

By the same argument, for $p_i\geq 2$ and $\delta_{i-1}= 0$ as well as when $p_i=1$
              \begin{equation}
            \label{20.104.1}
\frac{\eps^{-1}}{[\eps^{-2}(t_i-t_{i-1})]^{p_i/2}+1}
\leq
\frac{\eps^{-1}}{[\eps^{-2}(t_i-t_{i-1})]^{1/2-\zeta}}\leq
\eps^{-2\zeta}\frac{1}{(t_i-t_{i-1})^{1/2-\zeta}}
               \end{equation}
The product of all terms with powers of $\eps^{-2\zeta}$ is bounded by $\eps^{-2\zeta  M}$ uniformly in $\pi$.  The product of the $\eps$-factors gives:
             \begin{equation}
             \label{20.107}
C(\eps,\Delta):=
\eps^{|\{\delta_i=0,1\}|}
\Delta^{-|\{\delta_i=1\}|-\frac 12|\{\delta_i=0\}|}
\eps^{|\{\delta_i=2,\, p_i\geq 2,\delta_{i-1}\neq 0\}|}.
           \end{equation}
Hence, using arguments analogous the the ones used to get \eqref{20.104}, we have
              \begin{equation}
              \label{20.105}
\int_{0}^t dt_1\ldots\int_{t_{m-1}}^{t-\Delta} dt_m
\,f_{1,\ldots,m}
\leq c\,\eps^{-2\zeta  M}
C(\eps,\Delta)
\int_{0}^t dt_1\ldots\int_{t_{m-1}}^{t-\Delta} dt_m
\,\tilde f_{1,\ldots,m}
               \end{equation}
where the new integrand $\tilde f_{1,\ldots,m}$  is independent of $\eps$ and given by
              \begin{equation}
              \label{20.106}
\tilde f_{1,\ldots,m}:=
\prod_{\delta_i=1}\[
\frac{1}{(t_i-t_{i-1})^{1/2-\zeta}}
\]
\prod_{\delta_i=2}\[
\frac{1}{(t_i-t_{i-1})^{q_i-\zeta}}
\]
\prod_{\delta_i=0}\[
\frac{1}{(t_{i+1}-t_{i})^{1/2}}
\frac{1}{(t_i-t_{i-1})^{1/2-\zeta}}
\]
            \end{equation}
$q_i$ being defined for  $i:\delta_i=2$ as follows:  $q_i=\frac 12$ when either $p_i=1$ or $p_i\ge 2$ and $\delta_{i-1}=0$; in all other cases $q_i=1$, i.e.\ when  $p_i\ge 2$ and $\delta_{i-1}>0$.
To prove that the integral is finite we observe that for any $u<v$, $\alpha<1$, $\beta<1$,
            \begin{equation*}
\int_{u}^v   \frac 1{(s-u)^\alpha (v-s)^{\beta}} ds= c_{\alpha,\beta} (v-u)^{1-\alpha-\beta}
            \end{equation*}
with $\dis{ c_{\alpha,\beta}= \int_{0}^1   \frac 1{ s ^\alpha (1-s)^{\beta}}ds <\infty}$.  We use the above formula when integrating successively $t_m$, $t_{m-1},\ldots$ observing that the sum $\alpha+\beta$ at each step is $<1$ so that the resulting expression is bounded by a constant.  Once established that the integral is finite  a {\it scaling}
argument yields 
            \begin{equation}
            \label{20.108}
\int_{0}^t dt_1\ldots\int_{t_{m-1}}^t dt_m
\mathbf 1_{t_m<t-\Delta}
\,\tilde f_{1,\ldots,m}\leq c\,
t^{S}
            \end{equation}
where
            \begin{equation}
            \label{20.108.00}
S\ge  \frac 12 \Big( |\{\delta_i=1\}| +  |\{\delta_i=2, p_i=1\}| +  |\{\delta_i=2, p_i\ge 2,\delta_{i-1}=0\}|\Big) -\zeta  M
            \end{equation}
Proposition \ref{prop6.1} is proved.  \qed

\vskip2cm

   \section{Bounds when times cluster}
   \label{s5a}

The expectation which appears in  $w''_\pi(\und x,t)$ could be bounded exactly as in  $w'_\pi(\und x,t)$, the problem is that the power $S$ in \eqref{20.108} could then be negative and spoil the final bound.  For $w'_\pi(\und x,t)$
in fact the powers $(t-t_i)^{-1/2}$ and $(t-t_i)^{-1}$ could be bounded by $\Delta^{-1/2}$ and $\Delta^{-1}$ respectively, now $t-t_i$ might be smaller than $\Delta$. The factors  $(t-t_i)^{-1/2}$ and $(t-t_i)^{-1}$ may produce a negative $S$ in the integral in \eqref{20.108}.  The argument used in the proof  of Lemma \ref{lemma6e.1} to bound the expectation was based on an iterative argument where at each step we conditioned on the ``previous time'' bounding the conditional expectation uniformly on the positions of the particles at the conditioning time. We should here do better taking into account the fact that the conditioning configuration may be ``favorable''.

Using throughout the sequel the notation
$t_0\equiv 0$ and $t_{m+1}\equiv t$, we let, for $1\le H\le m$
            \begin{equation}
            \label{20.13}
\mathcal T_H:=\{
t_1,\ldots,t_m :\,t_{i+1}-t_i<\Delta,\,i=H,\ldots,m;\;\;
 t_H-t_{H-1}\ge \Delta
\}
             \end{equation}
and call  $t_H,\ldots,t_m,t$ the ``last cluster". Since we are supposing $t> (M+1)\Delta$, we readily see that the domain of integration $\{t_m\geq t-\Delta\}$ in $w''_\pi(\und x,t)$  can be decomposed as the union of $\mathcal T_H$ for $H=1,\dots m$.
In any $\mathcal T_H$ there is a time $t_{H-1}$ not belonging to the last cluster (this is $t_0\equiv 0$ if $H=1$).  Recalling that $A_i$ denotes the set of labels at time $t_i^+$ (the ``alive particles''),  given $H$ and  $t_1,\dots,t_m\in \mathcal T_H$,  we give the following definition.

\vskip.5cm

{\bf Definition.}  We denote by $w''_{\pi,H}(\und x,t)$ the integral in \eqref{20} extended to $t_1,\dots,t_m\in \mathcal T_H$. We then write
 $G_H$ for the set of all indices $i\ge H$ such that $\delta_i<2$ and there is $\ell\in A_{H-1} \cap\{k_i,l_i\}$ such that: $\ell \notin \{k_j,l_j\}$ for any $j\in [H,i)$ with $\delta_j=0$.

\vskip1cm

     \begin{prop}
     \label{prop7.1}
For any $\zeta>0$ there is $c$   so that for all $\pi:m=m(\pi) \le M$, for all
$\und x: |\und x|=n$, for all $\eps>0$ and all
$t \le \eps^{\beta^*}$, recalling \eqref{6e.0.4}, we have
          \begin{equation}
         \label{7e.0.3}
w''_{\pi,H}(\und x,t)
\leq \{c (\eps^2t)^{-\zeta  M} \Delta^{-S_1(H-1)} \eps^{S_2(H-1)} t^{S_3(H-1)}\}
\{c\, (\eps^2 \Delta)^{-\zeta  M} (\eps^{-2}\Delta)^{-\frac 14 |G_H|}
\Delta^{\frac 12|\{i\in [H, m]:\delta_i=2\}|}\}
       \end{equation}
 where the first curly bracket is equal to the right hand side of \eqref{6e.0.3} with $m$ replaced by $H-1$.

    \end{prop}

\vskip1cm

Recalling the definition of $\phi_h$ in \eqref{6e.1} we define for $h\ge H$
            \begin{equation}
         \label{6e.65}
 \phi^*_h:=
\prod_{i>h:\,i\in G_H}
\frac{1}{[\eps^{-2}\Delta]^{1/4-\zeta}+1}
\prod_{i>h:\,\delta_i=0}
\frac{1}{[\eps^{-2}(t_{i+1}-t_{i})]^{1/2}+1}
      \end{equation}
and set $\phi^*_h=\phi_h$ for $h<H$.  The analogue of Lemma \ref{lemma6e.1} then holds:

\vskip1cm
             \begin{lemma}
             \label{lemma6e.3}
There is $c$ and, for any $k$, $c'$ so that for all $H$, all $t_1,\dots,t_m\in \mathcal T_H$ and for any $h$ such that  $H\le h\le m$
               \begin{equation}
               \label{6e.7}
  \phi^*_h \EE_\eps \[ \psi_h\] \le c \phi^*_{h-1} \EE_\eps \[ \psi_{h-1}\] + c' (\eps^{-2}\Delta)^{-k}
       \end{equation}

               \end{lemma}

\medskip
{\em Proof.} Let  $h\ge H$. If $h\notin G_H$
and  $\delta_h<2$ we bound $\mathbf 1_{|x_{k_h}(t_h)-x_{l_h}(t_h)|=1}\le 1$. Analogously if  $\delta_h=2$ so that $h\notin G_H$,  we bound $\mathbf 1_{R_h}\le 1$ and in both cases we proceed as  in Lemma \ref{lemma6e.1}. It thus remains to consider the case $h\in G_H$ where the analysis is much more complex and relies heavily on the results in  Section \ref{secN6}.

Recalling \eqref{6e.1} and calling $\ell$ the particle-label entering in the definition of $G_H$, we  factorize $\psi_h  = \psi_{H-1} \psi^{(\ne \ell)}_{H,h} \psi^{( \ell,+)}_{H,h}$ where $\psi^{(\ne \ell)}_{H,h}$ does not depend on $x_\ell(\cdot)$ and it is given by
            \begin{equation}
                 \label{6e.8}
 \psi^{(\ne \ell)}_{H,h} :=
\prod_{H\le i< h:\delta_i=0,1}
\mathbf 1_{|x_{k_i}(t_i)-x_{l_i}(t_i)|=1}
\prod_{H\le i< h:\delta_i=0}\mathbf 1_{T_i}
\prod_{H\le i< h:\delta_i=2}  \mathbf 1_{ (x_j(t_i^-), j\ne \ell) \in I^c_{u_i}}
            \end{equation}
Then $\psi^{( \ell,+)}_{H,h}$ is equal to
            \begin{eqnarray}
 && \psi^{(\ell,+)}_{H,h} :=  \psi^{(\ell)}_{H,h} \mathbf 1_{T_h},\quad  \psi^{(\ell)}_{H,h}:=
 \mathbf 1_{|x_{k_h}(t_h)-x_{l_h}(t_h)|=1}
\prod_{ i\in \mathcal I}  \mathbf 1_{ x_\ell(t_i^-)\in I^c_{u_i}}\nn
\\
                \label{6e.9}
\\&&
\mathcal I =\big\{i\in [H,h]: \delta_i=2 \big\}\nn
            \end{eqnarray}
By \eqref{6e.2.0}--\eqref{6e.2.01} we get
               \begin{equation}
               \label{6e.10}
  \EE_\eps \[ \psi_h\] \le   \frac{c}{[\eps^{-2}(t_{h+1}-t_{h})]^{1/2}+1}\;  \EE_\eps \[ \psi_{H-1}  \psi^{(\ne \ell)}_{H,h}\psi^{(\ell)}_{H,h} \]
       \end{equation}
Since $h\in G_H$, $\ell$ is one of the two labels, $k_h,l_h$, for the sake of definiteness suppose
$l_h=\ell$. We compute the expectation on the right hand side of \eqref{6e.10} by conditioning on $\mathcal F(t_{H-1})$, then $ \psi_{H-1} $ drops from the conditional expectation and we are left with
                \begin{equation}
               \label{6e.11}
 \EE_\eps^* \big( \psi^{(\ne \ell)}_{H,h}\;\psi^{(\ell)}_{H,h} \big):=    \EE_\eps \Big[ \psi^{(\ne \ell)}_{H,h}\;\psi^{(\ell)}_{H,h}\;\big|\; \mathcal F(t_{H-1}) \Big]
       \end{equation}
Denote by $\und x^*$ the configuration at time $t_{H-1}^+$, $A_{H-1}$ being their labels; one of the particles is $\ell$ with position  $x^*_\ell$.  We shall estimate \eqref{6e.11}  using the coupling with independent particles defined in Section \ref{secN6}. More precisely in each interval $(t_i,t_{i+1})$, $i\ge H-1$,
we realize (couple) the stirring process in terms of independent random walks $\und x^0(s)$ adding and removing particles at the end of each interval as specified by $\pi$, independent particles are added on the same sites as the stirring ones, the sites are specified by the skeleton $\pi$; the independent particles start, as the stirring ones, at $t_{H-1}$  from $\und x^*$. The coupling is defined in terms of a priority list which is arbitrary except for the requirement that particle $\ell$ has the lowest priority.  The reason for this is that any stirring particles position $x_m(s)$, $m\ne \ell$,  $t_{H-1}\le s\le t$, is then a function of only  the independent processes $x^0_k(s)$, $t_{H-1}\le s\le t$, and $k\ne \ell$.
Recall that the stirring particles $x_j(s)$ are functions of the $\und x^0(\cdot)$, as described above, so that in the end the expression \eqref{6e.11}   becomes an expectation for independent particles.

We decompose the identity by writing $1= \chi + (1-\chi)$ where
                \begin{equation}
               \label{6e.12}
     \chi = \mathbf 1_{|x_\ell(t_h)-x^0_{\ell}(t_h)| \le (\eps^{-2}\Delta)^{1/4+\zeta}} \prod_{ i \in \mathcal I} \mathbf 1_{|x_\ell(t_i)-x^0_{\ell}(t_i)| \le (\eps^{-2}\Delta)^{1/4+\zeta}}
       \end{equation}

By \eqref{19.e11}  for any $n$ there is $c$ so that
               \begin{equation}
               \label{6e.13}
      \EE^*_\eps \Big[ \psi^{(\ne \ell)}_{H,h}\;\psi^{(\ell)}_{H,h}(1-\chi) \Big] \le (h-H+1) c (\eps^{-2}\Delta)^{-n}
       \end{equation}
Writing $F_i$ for all sites $x\in \La_N$ which have distance $>  (\eps^{-2}\Delta)^{1/4+\zeta}$ from $I_{u_i}$
we introduce the variable $\om$ with values in $\mathcal I$, see \eqref{6e.9}, which has value $i<h$ if $i$ is the smallest integer in $\mathcal I$ such that  $x^0_\ell(t_i) \notin F_i$. We set $\om=h$ if such $i$ does not exist.  We then have $\dis{1= \sum_{i\in \mathcal I} \mathbf 1_{\om = i}}$. We bound for $i<h$
               \begin{eqnarray}
               \label{6e.14}
       \EE^*_\eps \Big[ \psi^{(\ne \ell)}_{H,h} \psi^{(\ell)}_{H,h}\;\chi  \mathbf 1_{\om = i}\Big] &\le &
          \EE^*_\eps \Big[ \psi^{(\ne \ell)}_{H,h}\;    \mathbf 1_{x^0_{\ell}(t_i) \notin F_i}\Big]  =
              \EE^*_\eps \Big[ \psi^{(\ne \ell)}_{H,h}\Big]\mathcal P\big[ x^0_{\ell}(t_i) \notin F_i\big]
              \nn\\ &
              \le & \frac c{(\eps^{-2}\Delta)^{1/4-\zeta}+1} \EE_\eps \Big[ \psi^{(\ne \ell)}_{H,h}\;\big|\; \mathcal F(t_{H-1}) \Big]
                     \end{eqnarray}
When $\om =h$ we have for all $i\in \mathcal I, i<h$
              \begin{eqnarray}
               \label{6e.15}
    \mathbf 1_{x_\ell (t_i)\in I_{u_i^c} }\chi  \mathbf 1_{\om =h} = \chi  \mathbf 1_{\om =h}
       \end{eqnarray}
 and
               \begin{eqnarray}
               \label{6e.16}
    \mathbf 1_{|x_\ell(t_h)-x_{k_h}(t_h)|=1}  \chi  \mathbf 1_{\om =h} \le \chi  \mathbf 1_{\om =h}
     \mathbf 1_{|x^0_\ell(t_h)-x_{k_h}(t_h)|\le 1+ (\eps^{-2}\Delta)^{1/4+\zeta}}
       \end{eqnarray}
so that $\dis{\EE^*_\eps \Big[ \psi^{(\ne \ell)}_{H,h} \psi^{(\ell)}_{H,h}\;\chi  \mathbf 1_{\om = h}\Big] \le
\EE^*_\eps \Big[ \psi^{(\ne \ell)}_{H,h}\;     \mathbf 1_{|x^0_\ell(t_h)-x_{k_h}(t_h)|\le 1+ (\eps^{-2}\Delta)^{1/4+\zeta}} \Big] }$.
We condition on $\mathcal F^{(\ell)}$, the $\si$-algebra generated by the variables $x^0_j$, $j\ne \ell$, (including their ``attempted jumps'', see Section \ref{secN6}) observing that under $\PP_\eps[\cdot |\mathcal F^{(\ell)}]$ the variable $x^0_\ell$ is a simple random walk. Then
              \begin{eqnarray}
               \label{6e.17}
       \EE^*_\eps \Big[ \psi^{(\ne \ell)}_{H,h} \psi^{(\ell)}_{H,h}\;\chi  \mathbf 1_{\om = h}\Big] &\le &
                        \EE^*_\eps \Big[ \psi^{(\ne \ell)}_{H,h}\Big(\PP_\eps\Big[ |x^0_\ell(t_h)-x_{k_h}(t_h)|\le 1+ (\eps^{-2}\Delta)^{1/4+\zeta} \big|\mathcal F^{(\ell)}\Big]\Big)\Big]
              \nn\\ &
              \le & \frac c{(\eps^{-2}\Delta)^{1/4-\zeta}+1}  \EE_\eps \Big[ \psi^{(\ne \ell)}_{H,h}\;\big|\; \mathcal F(t_{H-1}) \Big].
       \end{eqnarray}
From \eqref{6e.10}, \eqref{6e.11}, \eqref{6e.14}, \eqref{6e.17}  and observing that $\EE_\eps \big(\psi_{H-1} \psi^{(\ne \ell)}_{H,h}\big)=\EE_\eps \big(\psi_{h-1}\big)$, we get  \eqref{6e.7}.  \qed

\vskip1cm

{\bf Proof of Proposition \ref{prop7.1}.}
Using Lemma \ref{lemma6e.3} we get
             \begin{equation}
             \label{23.2}
w''_{\pi,H}(\und x,t) \le   \int_0^t dt_1\ldots\int_{t_{m-1}}^{t} dt_m
\mathbf 1_{\mathcal T_H}
\,
f_{1,\ldots,H-1}
\,
g_{H,\ldots,m}
                 \end{equation}
where  $\mathcal T_H$ is defined in \eqref{20.13}, $f_{1,\ldots,H-1}$ in \eqref{20.7} and
                \begin{eqnarray}
                \label{23.1}
g_{H,\ldots,m} & := &
(\eps^{-2}\Delta)^{(-\frac 14+\zeta) |G_H|}
\prod_{\delta_i=1}\[
\frac{\eps^{-2}}{[\eps^{-2}(t-t_i)]^{1-\zeta}+1}
\]
\prod_{\delta_i=2}\[
\frac{\eps^{-1}}{[\eps^{-2}(t_i-t_{i-1})]^{p_i/2}+1}\]
\nonumber\\
&&
\prod_{\delta_i=0}\[
\frac{\eps^{-2}}{[\eps^{-2}(t-t_i)]^{1/2-\zeta}+1}
\frac{1}{
[\eps^{-2}(t_{i+1}-t_i)]^{1/2}+1}
\].
               \end{eqnarray}

We   bound all factors in $g_{h,\ldots,m}$
with $\delta_i=2$ as in \eqref{20.104.1}, in all the others we drop the $+1$ addendum in the denominator so that
we have a product of pure powers.  The same argument used in the proof of Lemma \ref{lemma6e.1} shows that the integral is finite. After observing that $t_H \ge t -(m-H+1)\Delta$,   a scaling argument yields:
              \begin{equation}
              \label{30.25}
           \text{   {\rm l.h.s. of \eqref{23.1}}}
\leq c\, (\eps^2 \Delta)^{-\zeta  M}
(\eps^{-2}\Delta)^{-\frac 14 |G_H|}
\Delta^{\frac 12|\{\delta_i=2\}|}
\end{equation}
which gives the second curly bracket in \eqref{7e.0.3}.  We are left with
    \begin{equation}
            \nn 
   \int_0^t dt_1\ldots\int_{t_{H-2}}^{t} dt_{H-1}\mathbf 1_{t_{H-1}<t-\Delta}
f_{1,\ldots,H-1}
                 \end{equation}
which is the same as in Proposition \ref{prop6.1}, thus giving the first curly bracket in \eqref{7e.0.3}. \qed

\vskip2cm

\section{Proof of Theorem~\ref{t1}}
  \label{s6}

Since $M$ is fixed the number of skeletons $\pi$ is finite. Hence recalling  \eqref{19}, \eqref{6e.0.2} and the definition of $w''_{\pi,H}$ before Proposition \ref{prop7.1}, it will suffice to prove
             \begin{equation}
              \label{8e.1}
\max_{\pi}w'_{\pi}(\und x,t)\le c (\eps^{-2}t)^{-c^* n},\quad  \max_{\pi,H}w''_{\pi,H}(\und x,t) \le c (\eps^{-2}t)^{-c^* n}.
          \end{equation}
We will need to distinguish various cases: firstly if $m(\pi)<M$ or $m(\pi)=M$; then if $t>(M+1)\eps^a$ (see \eqref{6e.0.0}) or the opposite and in all these cases we will have different arguments for $w'_{\pi}$ and $w''_{\pi,H}$.

\vskip1cm

$\bullet$\;{\bf {$m(\pi)=M$. Estimate of $w'_{\pi}(\und x,t)$ for $t>(M+1)\eps^a$.}}

From \eqref{6e.0.3}--\eqref{6e.0.4} we get
     \begin{equation}
              \label{a.9.2.0}
w'_{\pi}(\und x,t) \le c(\eps^2t)^{-\zeta M} [\eps\Delta^{-1}]^{|\delta_i=1|}
 \Delta^{-\frac12(|\delta_i=0|}\eps^{(|\delta_i=0|+|\delta_i=2,p_i\ge 2,\delta_{i-1}\ne 0|)} \, \,\,t^{\frac 12 |\delta_i=2, p_i=1|}
            \end{equation}
having used $t^{[S_3(m)}\le t^{\frac 12|\delta_i=2, p_i=1|]}$.
Observing that
             \begin{equation}
              \label{8e.1.0}
|\delta_i=0|+|\delta_i=2,p_i\ge 2,\delta_{i-1}> 0| \ge |\delta_i=2,p_i\ge 2|
            \end{equation}
we bound
      \begin{eqnarray*}
&&\hskip-1cm \Delta^{-\frac12(|\delta_i=0|} \eps^{(|\delta_i=0|+|\delta_i=2,p_i\ge 2,\delta_{i-1}\ne 0|)}\le \Delta^{-\frac12(|\delta_i=0|}
  \eps^{(|\delta_i=0|+\frac 12|\delta_i=2,p_i\ge 2,\delta_{i-1}\ne 0|} \eps^{\frac 12|\delta_i=2,p_i\ge 2,\delta_{i-1}\ne 0|}
 \\&&\le
 [ \eps\Delta^{-\frac12}]^{\frac12(|\delta_i=0|+|\delta_i=2,p_i\ge 2|)}
                 \end{eqnarray*}
We thus get
       \begin{equation}
             \label{8a.2.0}
w'_{\pi}(\und x,t) \le c(\eps^2t)^{-\zeta M} [\eps\Delta^{-1}]^{|\delta_i=1|} [\eps\Delta^{-\frac 12}]^{\frac 12(|\delta_i=0|+|\delta_i=2,p_i\ge 2|)} \, \,\,t^{\frac 12 |\delta_i=2, p_i=1|}
           \end{equation}

Let $\eps^b:= \max\{\eps\Delta^{-1};[\eps\Delta^{-\frac 12}]^{\frac12};\eps^{\frac{\beta^*}2}\}$ then
 \begin{equation}
             \label{8e.2.0}
w'_{\pi}(\und x,t) \le c(\eps^2t)^{-\zeta M} [\eps^b]^{|\delta_i=1|+|\delta_i=0|+|\delta_i=2,p_i\ge 2|+ |\delta_i=2, p_i=1|} \le (\eps^{2+\beta^*})^{-\zeta M} \eps^{bM}
           \end{equation}
We choose $M$ so that
             \begin{equation}
             \label{8e.2.1}
 bM \ge 2n
           \end{equation}
which implies
             \begin{equation}
             \label{8e.2.2}
w'_{\pi}(\und x,t) \le c(\eps^{2+\beta^*})^{-\zeta M} \eps^{2n} \le c' \eps^n
           \end{equation}
by choosing $\zeta$ small enough. Hence this term is compatible with \eqref{8e.1} provided $c^* < \frac 12 $.

\vskip1cm

$\bullet$\; {\bf  $m(\pi)=M$. Estimate of $w'_{\pi}(\und x,t)$ for $t\le (M+1)\eps^a$.}

In this case $t$ is proportional to $\Delta$ hence by \eqref{6e.0.3}--\eqref{6e.0.4}
             \begin{equation*}
w'_{\pi}(\und x,t) \le c(\eps^2t)^{-\zeta M} [\eps t^{-\frac 12}]^{|\delta_i=0|+|\delta_i=1|}\, \eps^{|\delta_i=2, p_i\ge 2,\delta_{i-1}>0|}\,\,t^{\frac 12 \big\{|\delta_i=2, p_i=1|+|\delta_i=2, p_i\ge 2,\delta_{i-1}=0|\big\}}
           \end{equation*}
We may suppose without loss of generality that $\eps^{-2}t\ge 1$, in such a case
             \begin{equation*}
w'_{\pi}(\und x,t) \le c\eps^{-4\zeta M}  \eps^{|\delta_i=2, p_i\ge 2,\delta_{i-1}>0|}\,\,t^{\frac 12 \big\{|\delta_i=2, p_i=1|+|\delta_i=2, p_i\ge 2,\delta_{i-1}=0|\big\}}
           \end{equation*}
           hence
                        \begin{equation*}
w'_{\pi}(\und x,t) \le c\eps^{-4\zeta M}  \eps^{\beta^*\big\{ |\delta_i=2, p_i\ge 2,\delta_{i-1}>0|+\frac 12 \big(|\delta_i=2, p_i=1|+|\delta_i=2, p_i\ge 2,\delta_{i-1}=0|\big)\big\}}
           \end{equation*}
             \begin{equation}
             \label{8e.2.3}
w'_{\pi}(\und x,t) \le c\eps^{-4\zeta M}  \eps^{\frac {\beta^*}2   |\delta_i=2|}
           \end{equation}

Observe that
             \begin{eqnarray}
            \nn
&&|\delta_i=1|+|\delta_i=0|+|\delta_i=2, p_i\ge 2|\le M
\\&& |\delta_i=1|+|\delta_i=0|\le n+(K-1)|\delta_i=2, p_i\ge 2|
         \label{20.107.02}
           \end{eqnarray}
so that
        \begin{equation}
             \label{20.107.031}
|\delta_i=1|+|\delta_i=0| \le n+\frac {K-1}K M
           \end{equation}
Since $|\delta_i=2|+|\delta_i=1|+|\delta_i=0|=M$, for $M> 2Kn$
       \begin{equation}
             \label{20.107.032}
|\delta_i=2|\ge M-\{  n+\frac {K-1}K M\} = \frac{M}{K}-n \ge \frac{M}{2K}
           \end{equation}
We choose $M$ so that
      \begin{equation}
             \label{20.107.03}
\frac {\beta^*}2 \frac{M}{2K} \ge 2n
           \end{equation}
Then
             \begin{equation}
             \label{8e.2.31}
w'_{\pi}(\und x,t) \le c\eps^{-4\zeta M}  \eps^{2n}
           \end{equation}
which for $\zeta$ small enough proves  compatibility with \eqref{8e.1} provided $c^*< \frac 12 $.

\vskip1cm

$\bullet$\;{\bf $m(\pi)=M$. Estimate of $w''_{\pi,H}(\und x,t)$ for $t>(M+1)\eps^a$.}

Suppose first $H\ge \frac M2$. We bound by 1 the second curly bracket in \eqref{7e.0.3} and use the same arguments
as those used to get \eqref{8a.2.0}. In this case, writing $|\delta_i=k|_H$  for $|\{i<H:\delta_i=k\}|$:
        \begin{equation*}
w''_{\pi,H}(\und x,t)
\leq  c (\eps^2t)^{-\zeta  M} [\eps \Delta^{-1}]^{|\delta_i=1|_H}[\eps\Delta^{-\frac 12}]^{\frac 12(|\delta_i=0|_H+|\delta_i=2,p_i\ge 2|_H)} \, \,\,t^{\frac 12 |\delta_i=2, p_i=1|_H}
       \end{equation*}
so that for $b$ as before \eqref{8a.2.0},
we have
          \begin{equation*}
w''_{\pi,H}(\und x,t)
\leq  c (\eps^2t)^{-\zeta  M}  \eps^{b(H-1)} \le (\eps^2t)^{-\zeta  M}  \eps^{b (M/2-1)}
       \end{equation*}
Taking $M$ such that $b (M/2-1)>2n$ we get, as in \eqref{8e.2.1},  $w''_{\pi,H}(\und x,t)\le c \eps^n$ compatible with \eqref{8e.1} provided $c^*< \frac 12 $.

  Suppose next that $H<\frac M2$.  We bound by 1 the first curly bracket in \eqref{7e.0.3} to write
          \begin{equation}
         \label{8e.2.4}
w''_{\pi,H}(\und x,t)
\leq  c\, (\eps^2 \Delta)^{-\zeta  M} (\eps^{-2}\Delta)^{-\frac 14 |G_H|}
\Delta^{\frac 12|\{i\in [H, m]:\delta_i=2\}|}
       \end{equation}

 We have that
      \begin{eqnarray}
           \nn
&&|G_H|\ge \frac 12 \big|\{i\ge H: \delta_i=1,0 \text{ and } \{k_i,\ell_i\}\cap A_{H-1}\ne \emptyset\}\big|
\\&& K|\{i\in [H,m]:\delta_i=2\}|\ge |\{i\in[H,m]:\delta_i=1,0 \text{ and } \{k_i,\ell_i\}\cap A_{H-1}= \emptyset\}|
  \label{8a.2.5}
           \end{eqnarray}
Leting $\eps^{b'}:= \max\{(\eps^{-2}\Delta)^{-\frac 18};\Delta^{\frac 12}\}$, one has
       \begin{equation}
         \label{8a.2.7}
w''_{\pi,H}(\und x,t)
\leq  c\, (\eps^2 \Delta)^{-\zeta  M} \eps^{b'[|i \in [H,m]: \delta_i=1,0 \text{ and } \{k_i,\ell_i\}\cap A_{H-1}\ne \emptyset|+|i \in [H,m]:\delta_i=2|]}
       \end{equation}
 First assume that (also simplifying a bit the notation)
       \begin{equation}
        \label{8a.2.7v}
|i \ge H: \delta_i=1,0 \text{ and } \{k_i,\ell_i\}\cap A_{H-1}= \emptyset|\ge \frac 12 |i \ge H: \delta_i=1,0|,
       \end{equation}
 and the second relation in \eqref{8a.2.5} yields
       \begin{equation*}
|i \ge H:\delta_i=2|\ge  \frac 1{2K} |i\ge H:\delta_i=1,0|
       \end{equation*}
 so that (since $4K >2$)
     \begin{eqnarray}
         \label{8a.2.7a}
&& |i \ge H:\delta_i=1,0 \text{ and } \{k_i,\ell_i\}\cap A_{H-1}\ne \emptyset|+|i \ge H:\delta_i=2|\ge |i \ge H:\delta_i=2|\nn\\&&
  \hskip2cm
  \ge
\frac 1{2} |i \ge H: \delta_i=2|+\frac 1{4K}  |i \ge H: \delta_i=1,0|\ge \frac {M-H}{4K}\ge
 \frac {M}{8K}
     \end{eqnarray}
From \eqref{8a.2.7} and if $ \frac {M}{8K}\ge 2n$, we get $|w''_{\pi,H}(\und x,t)|\le (\eps^2 \Delta)^{-\zeta  M} \eps^{b'2n}$,
compatible with \eqref{8e.1}.

For the case of complementary inequality in \eqref{8a.2.7v}, we have from \eqref{8a.2.7}:
       \begin{equation*}
w''_{\pi,H}(\und x,t)
\le  c\, (\eps^2 \Delta)^{-\zeta  M} \eps^{b'[\frac 12 |i \ge H: \delta_i=1,0|+|i \ge H: \delta_i=2|]} \le
c\, (\eps^2 \Delta)^{-\zeta  M} \eps^{b'[\frac 12 (M-H)]}
       \end{equation*}
with $\dis{M-H\ge \frac M2}$.

\vskip1cm

$\bullet$\; {\bf $m(\pi)=M$. Estimate of $w''_{\pi,H}(\und x,t)$ for $t\le (M+1)\eps^a$}.

 In  \eqref{7e.0.3} we estimate the first curly bracket as in  \eqref{8e.2.3} and we bound the second one with $ \Delta^{|i\ge H:\delta_i=2|}$, getting
             \begin{equation}
             \label{8e.2.8}
w''_{\pi,H}(\und x,t)
\leq  c\eps^{-4\zeta M}  \eps^{\frac {\beta^*}2   |i<H:\delta_i=2|} \Delta^{|i\ge H:\delta_i=2|}
           \end{equation}
and the same argument used for $m(\pi)=M$, $w'_{\pi}(\und x,t)$, $t\le (M+1)\eps^a$ applies.
We thus have compatibility with \eqref{8e.1} if $c^* < \frac 12 $.

\vskip1cm

$\bullet$\; {\bf $m(\pi)<M$. Estimate of $w'_{\pi}(\und x,t)$ for $t\le (M+1)\eps^a$}.

In this case $t$ is proportional to $\Delta$ hence by \eqref{6e.0.3}--\eqref{6e.0.4} we get after multiplying and dividing by $t^{-\frac 12 |\delta_i=2, p_i\ge 2,\delta_{i-1}>0|}$
             \begin{equation*}
w'_{\pi}(\und x,t) \le c(\eps^2t)^{-\zeta M} [\eps t^{-\frac 12}]^{|\delta_i=0|+|\delta_i=1|}\, [\eps t^{-\frac 12}]^{|\delta_i=2, p_i\ge 2,\delta_{i-1}> 0|}\,\,t^{\frac 12 \big\{|\delta_i=2, p_i=1|+|\delta_i=2, p_i\ge 2|\big\}}
           \end{equation*}
By using  the inequality
             \begin{equation}
              \label{8e.1.00}
n\le |\delta_i=0|+2|\delta_i=1| +K |\delta_i=2,p_i\ge 2|
          \end{equation}
and bounding $|\delta_i=2,p_i\ge 2|$ by  \eqref{8e.1.0}, we get
             \begin{equation}
             \label{8e.1.000}
n\le (K+1)|\delta_i=0|+2|\delta_i=1| +K |\delta_i=2,p_i\ge 2,\delta_{i-1}>0|
          \end{equation}
Notice that since $K\geq 1$ the largest factor is $K+1$.  Hence
           \begin{equation}
             \label{8e.2.10.00}
w'_{\pi}(\und x,t) \le c(\eps^2t)^{-\zeta M} [\eps t^{-\frac 12}]^{\frac{n}{1+K}}
\le c(\eps^{-2}t)^{-\frac{n}{2(1+K)}}
           \end{equation}
and for $\zeta$ small enough we have compatibility with \eqref{8e.1} if  $c^* < \frac{1}{2(1+K)} $.

\vskip1cm

$\bullet$\; {\bf $m(\pi)<M$. Estimate of $w'_{\pi}(\und x,t)$ for $t> (M+1)\eps^a$}.

By \eqref{6e.0.3}--\eqref{6e.0.4}, the inequality \eqref{8e.1.000}
and by dropping the factor  $t^{S_3(m)}$ we get
             \begin{eqnarray*}
w'_{\pi}(\und x,t) & \le &
c(\eps^2t)^{-\zeta M}
(\eps\Delta^{-1})^{|\delta_i=1|}(\eps\Delta^{-\frac 12})^{|\delta_i=0|}
\eps^{|\delta_i=2,\, p_i\ge 2,\,\delta_{i-1}> 0|}\\
&\le &
c(\eps^2t)^{-\zeta M}
\max\{(\eps\Delta^{-1})^{\frac 12},\,
(\eps\Delta^{-\frac 12})^{\frac{1}{K+1}},\,
\eps^{\frac 1K}
\}^n
\le
c(\eps^2t)^{-\zeta M} \eps^{\frac{1}{2(K+1)}n}
           \end{eqnarray*}
where the last inequality is true for our choice of $\Delta$ with $a=\frac{K}{K+1}\ge\frac{K-1}{K}$.
(Note that for $a\ge\frac{K-1}{K}$ the dominant term is $(\eps\Delta^{-1})^{\frac 12}$.)

Thus, for $\zeta$ small enough we have compatibility with \eqref{8e.1} if  $c^* < \frac{1}{4(K+1)} $.

\vskip1cm

$\bullet$\; {\bf $m(\pi)<M$. Estimate of $w''_{\pi,H}(\und x,t)$ for $t> (M+1)\eps^a$}.

First observe that since for all $\ell\in A_{H-1}$ there is $i\ge H$ such that $\delta_i=0,1$ and $ \ell\in \{k_i,\ell_i\}$ then
    \begin{equation}
    \label{9.25a}
|A_{H-1}|\leq 2 |G_H|+K |\delta_i=2|
        \end{equation}
Suppose first that $|A_{H-1}|<n$.  From \eqref{9.25a} we get that
the second factor in \eqref{7e.0.3} is bounded by
            \begin{equation}
             \label{8e.2.8.0}
c(\eps^2\Delta)^{-\zeta M}\(\max
\{
(\eps^{-2}\Delta)^{-\frac 18},\,\Delta^{\frac{1}{2K}}
\}
\)^{|A_{H-1}|}
\le
c(\eps^2\Delta)^{-\zeta M}
\eps^{\frac{1}{2(K+1)}|A_{H-1}|}
           \end{equation}
by our choice of $\Delta$, since for $a\le\frac{2K}{K+4}$ the dominant term is $\Delta^{\frac{1}{2K}}$.

On the other hand, for the first factor of \eqref{7e.0.3}
we apply the result as for the case
$m(\pi)<M$, $w'_{\pi}(\und x,t)$, $t> (M+1)\eps^a$ which gives the bound
$c(\eps^2t)^{-\zeta M} \eps^{\frac{1}{2(K+1)}(n-|A_{H-1}|)}$.

Overall we have:
             \begin{equation*}
w''_{\pi,H}(\und x,t) \le c(\eps^2\Delta)^{-2\zeta M}
\eps^{\frac{1}{2(K+1)}(n-|A_{H-1}|)}
\eps^{\frac{1}{2(K+1)}|A_{H-1}|}
              \end{equation*}
which for $\zeta$ small enough is  compatible with \eqref{8e.1} if  $c^* < \frac{1}{4(K+1)} $.
Same conclusions hold for $|A_{H-1}|\ge n$ as  it suffices to use only the bound \eqref{8e.2.8.0}.

\vskip1cm

$\bullet$\; {\bf $m(\pi)<M$. Estimate of $w''_{\pi,H}(\und x,t)$ for $t\le (M+1)\eps^a$}.
\nopagebreak

As in the previous case it suffices to investigate the case $|A_{H-1}|<n$.
The second factor in \eqref{7e.0.3} is  bounded by
           \begin{eqnarray*}
\le c(\eps^2\Delta)^{-\zeta M}\(\max
\{
(\eps^{-2}t)^{-\frac 18},\,t^{\frac{1}{2K}}
\}
\)^{|A_{H-1}|}
             \end{eqnarray*}
We have for $\dis{t\ge \eps^{\frac{2K}{K+4}}}$ (which is $\ge\eps^2$),
            \begin{eqnarray*}
\max
\{
(\eps^{-2}t)^{-\frac 18},\,t^{\frac{1}{2K}}
\}
 \le
t^{\frac{1}{2K}}
             \end{eqnarray*}
 For $\dis{\alpha:= \frac{1}{2(K+2)}}$

                         \begin{eqnarray*}
t^{\frac{1}{2K}} \le (\eps^2 t^{-1})^\alpha,\quad \text{for all}\;\; \eps^{\frac{2K}{K+4}}\le t\le \eps^a
             \end{eqnarray*}
%
%
%
For $\eps^{\frac{2K}{K+4}}\ge t \ge\eps^2$
             \begin{eqnarray*}
\max
\{
(\eps^{-2}t)^{-\frac 18},\,t^{\frac{1}{2K}}
\}
 \le
(\eps^{-2}t)^{-\frac 18},
             \end{eqnarray*}

By \eqref{8e.2.10.00} the first factor in \eqref{7e.0.3} is  bounded by
           \begin{equation}
             \label{8e.2.10.000}
 \le c(\eps^2t)^{-\zeta M} [\eps t^{-\frac 12}]^{\frac{n}{1+K}}
\le c(\eps^{-2}t)^{-\frac{n-|A_{H-1}|}{2(1+K)}}
           \end{equation}
Overall we have  compatibility with \eqref{8e.1} if  $c^* < \frac{1}{2(2+K)} $.  \qed

\vskip2cm

\section{Long times}
\label{sec:10}
\vskip.5cm

In this section we extend the estimate on the $v$-functions to times of order  $\log\eps^{-1}$. In Theorem \ref{t1} we have proved bounds for times $t\le \eps^{\beta^*}$, where $\beta^*$ is any given positive number. Here we shall study times $t\ge \eps^{\beta^*}$:

\vskip.5cm
    \begin{thm}
    \label{vteo}
For suitable $\tau>0$ (which depends on $c^*$ in Theorem \ref{t1})  and for any $n$ there is $c_n$ so that for all $\eps>0$
        \begin{equation}
        \label{it1}
\sup_{\eta_0}\,\sup_{\und x:\,|\und x|=n}\,\,\,\sup_{\eps^{\beta^*} \le t\le \tau\log\eps^{-1}}\,\,
|v^{\eps}(\und x,t|\eta_0)|\le c_n \eps^{(2-\beta^*)c^*\,n}
        \end{equation}
where  we have made explicit the dependence of the $v$-function on $\eps$ and on the initial configuration $\eta_0$.
    \end{thm}

\vskip.5cm

The theorem will be proved later in this section. We start with some definitions:

\vskip.5cm
 {\bf Definition}

  $\bullet$\; {\em
For any $b\in(0,1)$ and $f:\La_N\to \mathbb{R}$  we define
        \begin{equation}
        \label{i1}
|||f|||:=\sup_{x\in\La_N} |f|_x,\qquad |f|_x:=\big|\sum_{y\in\La_N} P^{(\eps)}_{\eps^{1+b}}(x, y)f(y)\big|
        \end{equation}
    $\bullet$\;
$\rho_\eps(x,t|f,s)$, $t\ge s$, denotes the solution of \eqref{pro.2.1} with initial datum $f$ at time $s$.}

\vskip.5cm

%

The following holds:
\vskip.5cm

            \begin{lemma}
            \label{lemma5.1}
 For any $\gamma>0$ such that $\dis{\ga\le \inf \big\{(2-\beta^*)c^*, \frac{1-b}4\big\}}$ the following holds.
 Given any $n$ there is $c_n$ so that for all $\eps$ and  for all configurations  $\eta$,
            \begin{equation}
        \label{i7}
\Pp_{\eps}^{\eta}\big(  |||\rho_\eps(\cdot, \eps^{\beta^*}|\eta,0)-\eta(\cdot,\eps^{\beta^*})|||\le \eps^{\ga}\big)\ge 1-c_n \eps^{n}
        \end{equation}
where $\eta(\cdot,\eps^{\beta^*})$ is the random configuration at time $\eps^{\beta^*}$ starting from $\eta$ at time $0$.
    \end{lemma}

    \vskip.5cm

{\bf Proof.} We first write
      \begin{equation*}
\Pp^\eta_{\eps}\big(  |||\rho(\cdot, \eps^{\beta^*}|\eta,0)-\eta(\cdot,\eps^{\beta^*})|||\ge \eps^{\ga}\big)\le
(2N+1)\sup_{x\in \La_N} \Pp^\eta_{\eps}\big( |\rho(\cdot,\eps^{\beta^*}|\eta,0)-\eta(\cdot,\eps^{\beta^*})|_x\ge \eps^{\ga}\big).
        \end{equation*}
By the Chebyshev inequality with power $2m$
     \begin{eqnarray*}
&&\hskip-2cm \Pp^\eta_{\eps}\big(  |\rho(\cdot,\eps^{\beta^*}|\eta,0)-\eta(\cdot,\eps^{\beta^*}) |_x\ge \eps^{\ga}\big)\le \eps^{-2m\ga }
 \sum_{x_1,., x_{2m}}\prod_{i=1}^{2m} P^{(\eps)}_{\eps^{1+b}}(x,x_i)\\&&\hskip7cm\E^{\eta}_\eps\big(
 \prod_{i=1}^{2m}[\eta(x_i,\eps^{\beta^*})-\rho(x_i,\eps^{\beta^*}|\eta,0)]\big).
    \end{eqnarray*}
There are constants $c$ and $c'$ (dependent on $m$) so that
     \begin{equation*}
|\E^{\eta}_{\eps}\big(
 \prod_{i=1}^{2m}[\eta(x_i,\eps^{\beta^*})-\rho(x_i,\eps^{\beta^*}|\eta,0)]\big)|\le c |v^\eps(Y,\eps^{\beta^*}|\eta)| \le c' \eps^{(-2+\beta^*)c^*|Y|}
    \end{equation*}
where $Y\subset\{x_1,\dots,x_{2m}\}$ is the set of singletons, i.e.\ all $x_i$ such that $x_j\ne x_i$ for $j\ne i$; the last inequality follows from  Theorem \ref{t1}.  On the other hand
     \begin{equation*}
P^{(\eps)}_{\eps^{1+b}}(x,x_i) \le c \big(\eps^{-2}\eps^{1+b}\big)^{-1/2} =c\eps^{(1-b)/2}
    \end{equation*}
and \eqref{i7} follows by the arbitrariness of $m$. \qed

\vskip1cm
We shall use the iterative approach used in  Chapter 5 of \cite{DP}, the presence of the factors $\eps^{-1}$ at the boundaries makes however the analysis different.
Thus, given $\tau>0$ (which will be specified later) we fix a time  $ \eps^{\beta^*}< T\le \tau\log\eps^{-1}$. Let $\dis{r\in[\frac 12 ,1]}$ be such that $T=(m+1)r\eps^{-\beta^*}, m\in \mathbb N$. We then partition the interval $\dis{[0,T]=\bigcup_{k=1}^{m+1} [t_{k-1},t_k)}$, $t_0=0$ where
 \begin{equation}
    \label{i4}
 [0,T]=\bigcup_{k=1}^{m+1} [t_{k-1},t_k),\;\;    t_k=kr\eps^{\beta^*},\quad m+1=r^{-1}\eps^{\beta^*}T. 
     \end{equation}
Given any initial configuration $\eta_0$, we denote by $\eta_k\equiv \eta(t_k)$ the random configurations at times $t_k$.

We study the evolution by successively conditioning the process at the times $t_k$. The conditioning at time $t_{k-1}$ fixes the configuration $\eta_{k-1}$ and by Lemma \ref{lemma5.1} the evolution in the next time step $[t_{k-1},t_k)$ is well approximated by $\rho_\eps(x,t|\eta_{k-1},t_{k-1})$.  To iterate this estimate we need to bound the difference
$\rho_{\eps}(x,t|\eta_{k}, t_k)-\rho_{\eps}(x,t|\eta_{k-1},t_{k-1})$: 
%
%

\vskip1cm
   \begin{prop}
    \label{rolemma}
Let $\gamma$ and $b$ be positive numbers as in Lemma \ref{lemma5.1}. Given any $\eta_{k-1}$, let  $\eta_k$ be such that $|||\rho_{k-1}(\cdot, t_k|\eta_{k-1})-\eta_k|||\le \eps^{\ga}$ then there is $c$ so that
     \begin{equation}
        \label{i9}
\big|\rho_\eps(x,t|\eta_k,t_k)-\rho_{\eps}(x,t|\eta_{k-1},t_{k-1})\big|\le c [\eps^\ga +\eps^b]\eps^{-\tau\pi}\tau\log\eps^{-1},\qquad  t\in(t_k, T]. 
        \end{equation}

    \end{prop}
\vskip.3cm

{\bf Proof.} We denote by
     \begin{equation}
 h_k(x,t)=|\rho_\eps(x,t|\eta_k,t_k)-\rho_{\eps}(x,t|\eta_{k-1},t_{k-1})|,\;\;\in \La_N;\quad h_k(t)=\sup_{x\in \La_N}h_k(x,t)
         \end{equation}
Since $\rho_\eps(\cdot)\le 1$, for any $t\in [t_k+\eps^{1+b},T]$ (using the weak form
of equation \eqref{pro.2.1})
    \begin{eqnarray}
    \label{r2}
&&h_k(x,t)\le \sum_y  P^{(\eps)}_{t-t_k-\eps^{1+b}}(x,y)|\rho_{\eps}(\cdot, t_k|\eta_{k-1},t_{k-1})-\eta_k|_y\nn\\&&\hskip2cm+c_1\eps^{-1}\int_{t_k}^{t}
\sum_{y\in I_+\cup I_-}  P^{(\eps)}_{s}(x, y)h_k(y,t-s)ds
    \end{eqnarray}
Since (recalling \eqref{a3.6})
    \begin{equation*}
\eps^{-1}\sum_{y\in I_+\cup I_-}  P^{(\eps)}_{s}(x, y)\le  \frac  {c_2 \tau \log\eps^{-1}}{\sqrt{s}}.
    \end{equation*}
We then bound
    \begin{eqnarray}
    \nn
&& \eps^{-1}\int_{t_k}^{t}
\sum_{y\in I_+\cup I_-}  P^{(\eps)}_s(x, y)h_k(x,t-s)\le 2K\eps^b+c_2 \tau \log\eps^{-1} \int_{t_k+\eps^{1+b}}^{t}
\frac 1{\sqrt s}h_k(t-s)ds.
    \end{eqnarray}
We thus get, using the hypothesis on $\eta_k$,
    \begin{equation}
    \label{r3}
h_k(x,t)\le [\eps^\ga +2K\eps^b]+c_3\tau \log\eps^{-1}\int_{t_k+\eps^{1+b}}^{t}
\frac 1{\sqrt s}h_k(t-s)ds\qquad t\in [t_k+\eps^{1+b},T].
    \end{equation}
Iterating \eqref{r3} and letting $a_n(t)$ be as in \eqref{a10.8}, we get that for all $t\in [t_k+\eps^{1+b},T]$
      \begin{eqnarray}
     \label{r4}
h_k(x,t)\le c_4[\eps^\ga +\eps^b](1+\tau \log\eps^{-1}\sum_{n=1}^\infty  a_n(t-t_k))
    \end{eqnarray}
Indeed the series converges as follows from Lemma \ref{sqrt} below. Since $t-t_k\le T$ $\le \tau\log\eps^{-1}$ from \eqref{a10.9a} we get
            $$
            \sum_{n=1}^\infty  a_n(t-t_k)\le
c e^{\pi\tau \log\eps^{-1}}= c\eps^{-\pi \tau}.
           $$
Thus from \eqref{r4} and \eqref{r6} we finally get \eqref{i9} for suitable $\tau$.
    \qed

  \vskip.5cm
 \begin{lemma}
    \label{sqrt}
    Let
        \begin{equation}
        \label{a10.8}
        a_n(t):=\int_0^{t}
\frac 1{\sqrt{s_1}}ds_1\int_0^{t-s_1}\frac 1{\sqrt{s_2}}ds_2\dots \int_0^{t-s_1\dots -s_{n-1}}
\frac 1{\sqrt{s_n}}ds_n.
    \end{equation}
Then
     \begin{equation}
        \label{a10.9}
a_n(t)\le (\pi t)^{\frac n 2} e^{-\frac n 2[\log (\frac n2)-1]}
    \end{equation}
and there is $c$ so that
    \begin{equation}
        \label{a10.9a}
 \sum_{n=1}^\infty
a_n(t)\le c \,  e^{\pi t}
    \end{equation}
  \end{lemma}
\vskip.3cm

{\bf Proof.} We have
     \begin{equation}
       \label{r5}
a_n(t)=\int_{[0,t]^n} \mathbf{1}_{s_1+\dots +s_n\le t}\prod_{i=1}^n \frac 1{\sqrt{s_i}}\,\,\,ds_1 \dots ds_n.
    \end{equation}
We change variable by setting $t_i=s_i\, t$ and get
     \begin{equation}
       \label{a10.10}
a_n(t)=(\sqrt t)^n \int_{[0,1]^n} \mathbf{1}_{t_1+\dots +t_n\le 1}\prod_{i=1}^n \frac 1{\sqrt{t_i}}\,\,\,dt_1 \dots dt_n.
    \end{equation}
Multiplying and dividing by $\exp \{-\alpha (t_1+\dots +t_n)\}$ we have
    \begin{equation}
       \label{r6}
a_n(t)\le (\sqrt t)^n e^{\alpha}\int_{[0,1]^n} \prod_{i=1}^n \frac {e^{- \alpha t_i}}{\sqrt{t_i}}\,dt_1 \dots dt_n\le
(\sqrt t)^n e^{\alpha}\Big[\int_0^1\frac {e^{- \alpha s}}{\sqrt{s}}\,ds\Big]^n\le
(\sqrt t)^n e^{\alpha}\big(\frac{\sqrt \pi}{\sqrt \alpha}\big)^n
    \end{equation}
    By choosing $\alpha=\frac n 2$ we get \eqref{a10.9} and \eqref{a10.9a} easily follows.\qed

 \vskip1cm

 {\bf Proof of Theorem \ref{vteo}}. Recalling the definitions given in \eqref{i4} and below and letting $\rho_0=\eta_0$, we write
    \begin{equation}
    \label{r8}
 -\rho_\eps(x,T|\eta_0,0)=\sum_{k=1}^m \big[ \rho_\eps(x,T|\eta_k,t_k)-\rho_{\eps}(x,T|\eta_{k-1},t_{k-1})\big]
 -\rho_\eps(x,T|\eta_m,t_m)
     \end{equation}
 Thus
            \begin{eqnarray*}
      \nn
&&v^\eps(X,T|\eta_0)=\E_{\eps}^{\eta_0}\Bigg(\prod_{x\in X}\Big[\big(\eta(x,T)-\rho_\eps(x,T|\eta_m,t_m)\big)\\&&\hskip4cm +\sum_{k=1}^m \big[ \rho_\eps(x,T|\eta_k,t_k)-\rho_{\eps}(x,T|\eta_{k-1},t{k-1})\big]\Big]\Bigg)
\\&&\hskip2cm =\sum_{Y\subset X}\E_{\eps}^{\eta_0}\Bigg(\prod_{y\in Y}\big(\eta(y,T)-\rho_\eps(y,T|\eta_m,t_m)\big)\\&&\hskip4cm \times
\prod_{x\in X\setminus Y}\sum_{k=1}^m \big[ \rho_\eps(x,T|\eta_k,t_k)-\rho_{\eps}(x,T|\eta_{k-1},t_{k-1})\big]\Big]\Bigg)
\nn
    \end{eqnarray*}
 We take conditional expectation with respect to $\mathcal F(\{\eta_t, t\le t_m\})$ and since the bound on the $v$ function in  Theorem \ref{t1} is uniform in the initial configuration ($\eta_m$ in this case) we get
     \begin{eqnarray}
      \nn
&&|v(X,t)|\le \sum_{Y\subset X}c_Y\,\eps^{c^*\,|Y|}\E_{\eps}^{\eta_0}\Bigg(\prod_{x\in X\setminus Y}\sum_{k=1}^m \big| \rho_\eps(x,T|\eta_k,t_k)-\rho_{\eps}(x,T|\eta_{k-1},t_{k-1})\big|\Bigg).
\nn
    \end{eqnarray}
 Let
    \begin{equation}
    \nn
 \mathcal{A}_\ga:=\{(\eta_1,\ldots,\eta_m): \sup_k |||\rho_{\eps}(\cdot, t_k|\eta_{k-1},t_{k-1})-\eta_k|||\le \eps^{\ga}\}.
    \end{equation}
 Since
    $$
  \Pp^{\eta_0}_\eps\Big(\sup_k |||\rho_{\eps}(\cdot, t_k|\eta_{k-1},t_{k-1})-\eta_k|||\ge \eps^{\ga}\Big)\le  m \sup_\eta
 \Pp^{\eta}_\eps\big(  |||\rho_\eps(\cdot, \eps^{\beta^*}|\eta,0)-\eta_{\eps^{\beta^*}}\||\le \eps^{\ga}\big)
    $$
and by Lemma \ref{lemma5.1} for any integer $q$ there is $c_q$ so that
    \begin{equation}
    \nn
  \Pp^{\eta}_\eps\big( \mathcal{A}_\ga\big)\ge 1-c_q \eps^q,
    \end{equation}
so that, using \eqref{i9} we get
             \begin{equation}
      \label{r9}
|v(X,t)|\le \sum_{Y\subset X}c_Y\,\eps^{c^*\,|Y|}\Big\{\bar c_q \eps^q
+[m\eps^{-d}(\eps^\ga+\eps^b)]^{X\setminus Y}\Big\}.
    \end{equation}
 We choose the parameters in such a way that
      \begin{equation}
      \label{r10}
\eps^{-\beta^*}d\log\eps^{-1}\eps^{-d}(\eps^\ga+\eps^b)\le \eps^{c^*}
    \end{equation}
 Theorem \ref{vteo} then follows from \eqref{r9} and \eqref{r10}. \qed

\vskip2cm

 {\bf Acknowledgments.}

 The research has been partially supported by PRIN 2007 (prot.20078XYHYS-003). M.E.V is partially supported by CNPq grant 302796/2002-9.
 M.E.V. thanks  Universit\`a di Roma ``Tor Vergata" and Universit\`a di Roma ``La Sapienza", and Universit\`a de
 L'Aquila for the support and hospitality during the visits when this work was carried out.
The research of D.T. has been supported by a Marie Curie Intra
European Fellowship within the 7th European Community Framework
Program.

\vskip1cm

\bibliographystyle{amsalpha}

\vskip1cm

\end{document}